\documentclass{article}

\usepackage{acronym}
\usepackage{amsmath}
\usepackage{amssymb}
\usepackage{booktabs}
\usepackage{graphicx}
\usepackage{subcaption}
\usepackage{siunitx}
\usepackage{fullpage}
\usepackage{amsthm}
\usepackage{authblk}
\usepackage{xr}
\usepackage{bm}
\usepackage{placeins}

\usepackage{booktabs,siunitx}

\theoremstyle{remark}
\newtheorem{remark}{Remark}

\newacro{ODIL}{optimizing a discrete loss}
\newacro{GN}{Gauss--Newton}
\newacro{LM}{Levenberg--Marquardt}
\newacro{CLM}{Constantin--Lax--Majda}
\newacro{CCF}{C\'{o}rdoba--C\'{o}rdoba--Fontelos}
\newacro{PINN}{physics-informed neural network}
\newacro{PDE}{partial differential equation}

\newcommand{\Hop}{\mathcal{H}}

\title{Discovery of Singularities in Partial Differential Equations\\ Using Adaptive Numerical Discretizations}
\date{}

\author[]{Allen Yang}
\author[]{Lucas Amoudruz}
\author[]{Petros Koumoutsakos}
\affil[]{Computational Science and Engineering Laboratory,\\ 
 Harvard University, 23 Oxford St, Cambridge, MA 02138, USA.}

\begin{document}
\maketitle

\begin{abstract}
We compute self-similar singularity profiles of \acp{PDE} using adaptive numerical discretizations based on \ac{ODIL}. The fields are represented in overlapping local Chebyshev patches joined with a mapped asymptotic patch. After factoring out the leading power laws of solution ansatzes, their limiting amplitudes at infinity are determined by the optimization, eliminating domain truncation and prescribed far-field conditions. Adaptive refinement is guided by Chebyshev-coefficient decay, and the nonlinear least-squares problem is solved using automatically differentiated Gauss--Newton iterations with Levenberg--Marquardt damping. The method achieves \(O(10^{-14})\)-level accuracy for Burgers profiles, residuals of \(O(10^{-14})\) to \(O(10^{-10})\) for the generalized De Gregorio family, and of \(O(10^{-8})\) in a local region for the two-dimensional Boussinesq equations. These results demonstrate that \ac{ODIL} can identify singularity profiles for \acp{PDE} with high accuracy using adaptive numerical methods, complementing ongoing efforts that rely on artificial intelligence. 

\end{abstract}
\acresetall

\section{Introduction}\label{sec:intro}
Finite-time singularities arise in many nonlinear \acp{PDE}, where they appear as unbounded amplitudes or gradients, or as the collapse of a spatial scale. Although their physical origins differ, such processes typically exhibit a self-similar structure near the singularity~\cite{eggers2008role,luo2014potentially}. Examples include blowup in semilinear heat equations~\cite{collot2019stability}, collapse in the focusing nonlinear Schr\"odinger and complex Ginzburg--Landau equations~\cite{dahne2024self}, and rupture in thin-film models~\cite{dallaston2017self}. In fluid mechanics, numerical studies have been central to identifying candidate singularity scenarios, and computer-assisted proofs have since turned high-accuracy numerical profiles into rigorous blowup results for the two-dimensional Boussinesq equations and the three-dimensional axisymmetric Euler equations with boundaries~\cite{luo2014toward,chen2022stable,chen2025stable}.

In recent years, inverse methods such as \acp{PINN}~\cite{Lagaris1997,raissi2019physics} have been used extensively in the discovery of smooth self-similar profiles for the two-dimensional Boussinesq equations, three-dimensional axisymmetric Euler equations with a boundary, the incompressible porous media equations, and the \ac{CCF} model~\cite{wang2023asymptotic,wang2025discovery,wang2025high}. Because these methods solve the stationary profile equations directly, they can also recover unstable branches that are difficult to obtain through forward time evolution. More recent work has produced a candidate self-similar profile for the unforced Euler equations on $\mathbb{R}^3$, together with a framework for proving its nonlinear stability, as part of an effort to resolve the longstanding open problem of singularity formation from smooth initial data with no forcing~\cite{ganeshram2026stable,duruisseaux2026stability}. These advances demonstrate the effectiveness of residual-based optimization for singularity discovery. 

In this paper, we demonstrate that residual-based optimization for singularity formation can be studied using classical numerical approximations whose resolution and differentiation properties are explicitly controlled.
We develop an adaptive spectral formulation of \ac{ODIL}~\cite{karnakov2024solving,karnakov2023flow} for computing self-similar singularity profiles. We represent unknown fields using polynomial interpolation so that the optimization variables are the nodal field values. For a fixed discretization, field values and spatial derivatives are linear functions of the optimization variables. In contrast to \acp{PINN}, the proposed parameterization does not introduce additional nonlinearity beyond that of the profile equations, while retaining the flexibility of \acp{PINN} in terms of residual-based optimization.

This direct parameterization also offers  advantages for constructing singularity profiles and assessing them through systematic numerical analysis. Approximate solutions immediately provide initial guesses through their nodal values, and interpolation transfers computed profiles between resolutions, supporting systematic refinement and continuation. Prescribed nodal values can be imposed by eliminating variables, while derivative and integral constraints can be expressed directly in terms of the discrete fields. Furthermore, the number of unknowns is determined by the required spatial resolution and can be chosen separately from the number of residual evaluation points. Enforcing the equations on a denser grid produces an overdetermined nonlinear least-squares problem without enlarging the representation itself. The Chebyshev-coefficient decay provides  diagnostics of spatial resolution, complementing independent residual evaluations and comparisons between successively refined solutions. Together, these features provide a framework for examining how the computed profiles, their derivatives, and their asymptotic behavior depend on numerical resolution---an essential consideration when using approximate profiles to investigate singularity formation. Furthermore, the small number of unknowns in the spectral representation substantially reduces the training time relative to \acp{PINN}. This enables computationally efficient studies for continuation and sensitivity. We assess the method on the one-dimensional Burgers and  generalized De Gregorio equations, and in two-dimensional Boussinesq profiles. These problems respectively provide an exact-solution benchmark, a nonlocal Hilbert-transform coupling, and a multidimensional system with coupled fields and angular far-field structure. 

The remainder of this paper is organized as follows. Section~\ref{sec:methodology-general} presents the methodology: Chebyshev interpolation, the multi-patch discretization in one dimension and its extension to two dimensions, the training loss, and the optimization. Section~\ref{sec:general-adaptive} describes the adaptive refinement of the patch decomposition. Section~\ref{sec:results} reports results for the Burgers equation (Section~\ref{sec:burgers}), the generalized De Gregorio family (Section~\ref{sec:degregorio}), and the two-dimensional Boussinesq equations (Section~\ref{sec:boussinesq}), and Section~\ref{sec:conclusions} concludes. The supplementary materials contain the derivation of the mapping exponent, the parity-extension procedure, implementation and configuration details for each application, and preliminary results for the three-dimensional Euler equations.

\section{Methodology}\label{sec:methodology-general}

In self-similar variables, each problem considered here reduces to a system
of nonlinear profile equations for stationary fields on an unbounded domain.
Symmetries of the profiles allow us to restrict the computational domain to
the half-line $[0,\infty)$ for the Burgers and generalized De Gregorio
problems, and to the quarter plane $[0,\infty)^2$ for the Boussinesq problem,
with the full profiles recovered by parity extension. The profiles are
smooth, may contain sharp features localized in small regions of the domain,
and approach known power laws in the far field.

We compute these profiles with \ac{ODIL}~\cite{karnakov2024solving}, which
minimizes the residuals of the governing equations directly over the values
of a discretized solution. Each unknown field is represented by overlapping
local Chebyshev patches together with a dedicated asymptotic patch for the
far field, blended by a smooth partition of unity. The nodal values on all
patches serve as the optimization variables, the residuals of the profile
equations are evaluated at collocation points on each patch, and the
resulting nonlinear least-squares problem is solved with damped
Gauss--Newton iterations. The patch decomposition and polynomial degrees can
then be adapted using the decay of the Chebyshev coefficients.

We first develop the method on the half-line.
Section~\ref{sec:cheb-background} reviews Chebyshev interpolation, introduces
the multi-patch representation, and constructs the residual-based loss and
its optimization on a fixed discretization.
Section~\ref{sec:methodology-general-2d} extends the representation to the
quarter plane, and Section~\ref{sec:general-adaptive} describes how the
discretization is adapted during optimization. These sections present the
main ideas and the details needed to follow the results, while the full
description of each component is given in the supplementary materials.
Equation-specific material, including the self-similar ansatz and constraints
of each problem and the evaluation of the Hilbert transform in the De Gregorio
residual, is given in Section~\ref{sec:results} and the supplementary
materials.

\subsection{Spectral formulation of \ac{ODIL} using Chebyshev polynomial interpolation}\label{sec:cheb-background}

The steady-state self-similar profiles considered in this paper are smooth,
and we exploit this smoothness through high-order polynomial interpolation.
For sufficiently smooth functions and a well chosen polynomial basis, the
interpolation error decays super-algebraically with the number of
interpolation points \cite{trefethen2013approximation,boyd2001chebyshev}, so
such profiles can be resolved with far fewer unknowns than a finite difference
discretization or a \ac{PINN}-based approach would require. We adopt Chebyshev
polynomials, which are well suited to the profiles studied here, although the
framework is agnostic to the choice of basis and could accommodate any suitable alternative.

We first review the properties of Chebyshev interpolation that are relevant for \ac{ODIL}. For a smooth scalar field
$q\in C^\infty([-1,1])$ and $n\geq 2$, let $q_n$ denote the unique
degree-$(n-1)$ polynomial interpolating $q$ on the $n$-point Chebyshev--Lobatto
grid (eq.~\eqref{sm-eq:cheb-nodes}). Because these points cluster near both
endpoints, $q_n$ avoids the Runge phenomenon and converges to $q$
super-algebraically as $n$ increases
\cite[Thm.~7.2]{trefethen2013approximation}. Equivalently, $q_n$ can be written
as a linear combination of the first $n$ Chebyshev basis polynomials:
\begin{equation}
  q_n(y) = \sum_{k=0}^{n-1} a_k T_k(y), \qquad T_k(y) = \cos(k\arccos y).
  \label{eq:cheb-expansion}
\end{equation}
This form is useful for both computation and error assessment: a fast discrete
cosine transform (an FFT-based algorithm) converts between value space (the set of
values of $q_n$ on the interpolation nodes) and coefficient space (the set of $a_k$) in
$\mathcal{O}(n\log n)$ operations in either direction; derivatives of $q_n$ can be
obtained exactly through a recurrence in coefficient space. The   magnitude of the trailing
coefficients $|a_k|$ provides an estimate of the interpolation error. Furthermore,
because $q_n$ is a polynomial, it is defined on the entire interval, not only at
the interpolation nodes, removing  a limitation of finite-difference
\ac{ODIL}, in which evaluating the solution off-grid requires a secondary interpolation step. Given the values of $q_n$ at the interpolation nodes, we use 
barycentric interpolation~\cite{berrut2004barycentric}, 
which is numerically stable even for large $n$ and does not require
forming the Chebyshev coefficients explicitly.
For Chebyshev--Lobatto points, the barycentric weights are known in closed
form, and when the evaluation points are also fixed, the full evaluation
matrix can be precomputed and reused, reducing the cost of the required repeated residual
evaluations.

We note that a single Chebyshev interpolant on $[-1,1]$ is insufficient for the problems considered in this paper. First, the self-similar profiles
are posed on the semi-infinite domain $[0,\infty)$ and may decay slowly or
even grow algebraically as $y\to\infty$; thus, the domain cannot simply be
truncated to a bounded interval. Second, when profiles contain sharp features
localized in small regions of the domain, a single interpolant is unable to resolve them 
efficiently. A piecewise representation can concentrate resolution where it is needed, but
when seeking globally smooth self-similar blowup profiles, the ansatz should be globally smooth as well. We therefore use overlapping patches
combined with a smooth partition of unity, along with a dedicated
asymptotic patch for the far field. This yields a composite representation
that resolves localized features efficiently, captures the far-field
behavior, and remains smooth across patch boundaries (full details are given in
Supplementary Section~\ref{sm-sec:grid}).

Each scalar field $q\colon[0,\infty)\to\mathbb{R}$ is represented on an
overlapping decomposition of the domain. The finite region is covered by $m$
local patches $K_i = (y_{\mathrm{start}}^{(i)}, y_{\mathrm{end}}^{(i)})$,
$i = 1,\dots,m$, which overlap their neighbors and together cover
$(0, y_{\mathrm{end}}^{(m)})$. On each $K_i$, $q$ is approximated by the
Chebyshev interpolant $q_{n_i}^{(i)}$ on an $n_i$-point Chebyshev--Lobatto grid
mapped to $\overline{K_i}$. The far field is covered by a single asymptotic
patch $K_{\mathrm{asym}} = (y_{\mathrm{start}}^{\mathrm{asym}}, \infty)$, with
$y_{\mathrm{start}}^{\mathrm{asym}} < y_{\mathrm{end}}^{(m)}$ so that it
overlaps the last local patch. Following~\cite{wang2025discovery}, on this
patch we factor out the leading-order behavior $E_q(y) = |y|^\sigma$, with
$\sigma$ determined by the governing equations, and approximate the normalized
field
\begin{equation}\label{eq:qhat}
  \widehat q(y) = \frac{q(y)}{E_q(y)},
\end{equation}
which approaches a constant even when $q$ grows at
infinity. To approximate this normalized field on the semi-infinite domain, we introduce 
a reference variable $x \in [-1,1]$ through an
algebraic map~\cite{boyd2001chebyshev} from the closure of the asymptotic patch $K_{\mathrm{asym}}$ onto $[-1,1]$, which sends $y=\infty$ to $x=1$, and approximate $\widehat q(y(x))$ by the Chebyshev
interpolant $\widehat q_{n_{\mathrm{asym}}}$ on an $n_{\mathrm{asym}}$-point
Chebyshev--Lobatto grid in $x$. The map depends on a mapping exponent
$\alpha>0$, which is chosen so that the normalized fields are smooth in $x$,
preserving the rapid convergence of the interpolant (Supplementary
Sections~\ref{sm-sec:grid} and~\ref{sm-sec:exponents}). The map also sends a
prescribed reference point $y_{\mathrm{mid}}^{\mathrm{asym}}$ to $x=0$, so that
half of the nodes of the asymptotic grid lie in
$[y_{\mathrm{start}}^{\mathrm{asym}}, y_{\mathrm{mid}}^{\mathrm{asym}}]$ and half
beyond it.

The patchwise interpolants are blended with a smooth partition of unity
$\{\chi^{(i)}\}_{i=1}^m \cup \{\chi^{(\mathrm{asym})}\}$ subordinate to the
patch covering, giving the composite approximation
\begin{equation}\label{eq:pou-blend}
  q_\mathbf{n}(y)
  =
  \sum_{i=1}^{m}\chi^{(i)}(y)\,q_{n_i}^{(i)}(y)
  +
  \chi^{(\mathrm{asym})}(y)\,\widehat{q}_{n_{\mathrm{asym}}}(y)\,E_q(y).
\end{equation}
Since each interpolant and each window function is smooth, so is
$q_\mathbf{n}$. The optimization variables are the nodal values on all
patches. Derivatives of $q_\mathbf{n}$ follow from Chebyshev differentiation of
the patchwise interpolants combined with the chain and product rules, so the
fields and their derivatives can be evaluated exactly (up to rounding) at any
point on the half-line $[0, \infty)$.

With this representation in place we now formulate the spectral \ac{ODIL} framework  (full details in Supplementary Section~\ref{sm-sec:loss}).
The unknowns are the nodal values of each field
on the local patches and of its normalized counterpart $\widehat q$ on the
asymptotic patch. Denoting by $f[\,\cdot\,](y)$ the residual of a governing
equation, we evaluate $f$ on each local patch $K_i$ at $\kappa n_i$
Chebyshev--Lobatto points, where the oversampling factor $\kappa \ge 1$ (we
use $\kappa = 2$ or $3$) improves residual control between nodes. On the
local patches, we also penalize the derivative of the residual, which we found
facilitates convergence to the desired profiles. On the asymptotic patch, we
evaluate $f$ at the finite nodes of the grid, excluding the node at
$y = \infty$. Finally, we penalize disagreement between neighboring
interpolants on their overlaps. Without these matching terms, compensating
changes to neighboring interpolants could leave the blended field nearly
unchanged, which creates near-null directions in the optimization. Individual
problems may add constraints or adjust the residual weights, as specified in
the corresponding sections.

The result is a nonlinear least-squares problem for  the nodal values and any additional unknown  parameters of the \acp{PDE}. The spectral representation keeps the
number of unknowns small and we can use \ac{GN} iterations with
Levenberg--Marquardt (\ac{LM})
damping~\cite{levenberg1944method,marquardt1963algorithm}, implemented in
\textsc{Optimistix}~\cite{rader2024optimistix} and
\textsc{Lineax}~\cite{rader2023lineax} within \textsc{JAX}~\cite{jax2018github}
in double precision. In our experiments, first-order optimizers did not
reliably reach the desired accuracy, while the default \ac{LM} settings
required almost no problem-specific tuning. Each damped subproblem is solved
by dense QR factorization of the residual Jacobian, except at the finest
Boussinesq resolution, where we use LSMR~\cite{fong2011lsmr} with randomized
preconditioning~\cite{avron2010blendenpik} to avoid storing the full
Jacobian.

\subsection{Extension to two dimensions}\label{sec:methodology-general-2d}

For the two-dimensional problem, we extend this representation to fields on
the quarter plane $[0,\infty)^2$ (full details in Supplementary
Section~\ref{sm-sec:2d}). Near the origin, we use Cartesian coordinates
$(y_1, y_2)$ and cover the region with a chain of rectangular patches that
overlap their neighbors in $y_1$ and share a common extent in $y_2$. On each
patch, $q$ is approximated by a tensor-product Chebyshev interpolant. This
arrangement is well suited to  the problems considered herein, but more general coverings
with patches overlapping in both directions can be handled similarly, as
in~\cite{aiton2019adaptive}. In the far field, mapping each Cartesian coordinate to $[-1,1]$
separately can reduce regularity at the corner where both $y_1$ and $y_2$
tend to infinity. We therefore represent the far field in polar coordinates
$(r,\theta)$, with one unbounded direction. The far-field region
$r > r_{\mathrm{start}}$ is split in $\theta$ into overlapping angular
patches. On each patch, we factor out the leading-order behavior $E_q(r) = r^\sigma$
and approximate the normalized field $\widehat q = q/E_q$ by a tensor-product
Chebyshev interpolant in $(\rho, \theta)$, where $\rho \in [-1,1]$ is obtained
from $r$ through the same algebraic map as in Section~\ref{sec:cheb-background}.

The result is two families of patches, a chain in $y_1$ near the origin and a
chain in $\theta$ in the far field, which we combine in two steps. First,
within each family, the patch interpolants are blended using 
one-dimensional partition-of-unity functions applied in $y_1$ for the local patches and in
$\theta$ for the asymptotic patches. Second, the two families are combined
with a radial partition of unity that transitions between $r_{\mathrm{start}}$
and a prescribed radius $r_{\mathrm{blend}}$, over which every point belongs
to both regions. The resulting composite approximation is smooth, and the
remaining components of the method, including the choice of unknowns, the
residual blocks, and the \ac{GN}/\ac{LM} optimization, carry over directly
from the one-dimensional case.

\subsection{Adaptivity}\label{sec:general-adaptive}

The representations developed in the preceding sections leave open how the
discretization itself should be chosen: the patch decomposition, the
overlaps, and the number of nodes on each patch must all be specified, and a
poor choice can either leave sharp features under-resolved or waste unknowns
where the profile varies slowly. The Chebyshev representation offers solutions to these challenges. The decay of the
Chebyshev coefficients on each patch indicates whether that patch is
adequately resolved, and this information can guide changes to the patch
decomposition and polynomial degrees. We remark that such a direct connection
between resolution assessment and refinement is not readily available for a standard neural-network parameterization.

Following the split, merge,
and chop approach in~\cite{aiton2017adaptive}, we subdivide patches where
additional resolution is needed, combine neighboring patches
where a single interpolant suffices, and reduce polynomial
degrees where fewer coefficients are sufficient. By
combining these operations with the \ac{ODIL} optimization,
we adapt the discretization as the approximate solution
develops. In this paper, we apply this procedure only to the local
patches of the one-dimensional problems, keeping the width
of the local region (i.e., $y^{(m)}_\mathrm{end}$; see Section~\ref{sec:cheb-background}) and the asymptotic patch's mapping and
resolution fixed. The local patches retain their
ordered, nearest-neighbor overlap structure throughout,
allowing us to manage them as a simple list rather than
a binary tree. This adaptivity allows us to efficiently resolve difficult features in regions of interest.

The resolution criterion that guides the adaptive refinement employed in this paper
is motivated by the relationship
between Chebyshev-coefficient decay and interpolation error.
In practice, we examine the coefficients
$a_0,\dots,a_{n-1}$ of each computed local interpolant and
define the normalized tail maximum
\[
  e_k
  =
  \frac{\max_{k\leq j\leq n-1}|a_j|}
       {\max_{0\leq j\leq n-1}|a_j|},
  \qquad k=0,\dots,n-1,
\]
taking $e_k=0$ when all coefficients vanish. For a prescribed
relative tolerance $\varepsilon$, we seek the smallest
admissible cutoff index $k$ satisfying $e_k\leq\varepsilon$.
We then retain $a_0,\dots,a_{k-1}$, subject to a prescribed
minimum resolution. This tail-based rule is a simplified
alternative to Chebfun's
\texttt{standardChop}~\cite{aurentz2017chopping}.

We use this coefficient-based diagnostic to guide both
patch refinement and coarsening, alternating optimization
with split, merge, and chop passes that alter the discretization.
Given a tolerance $\varepsilon$, a maximum polynomial degree
$N_{\max}$, and an overlap fraction $\eta\in(0,1)$, an iteration of the adaptive
algorithm proceeds as follows:
\begin{enumerate}
  \item \emph{Optimize.}
  Using~\ac{ODIL} as described in Section~\ref{sec:methodology-general},
  solve the problem of interest on the current discretization, taking the
  previously converged solution as the initial guess when available.
  To assess resolution on each patch $K_i$, sample the composite
  field~\eqref{eq:pou-blend} on the $N_{\max}$-node
  Chebyshev--Lobatto grid associated with $\overline{K_i}$ and
  compute the corresponding Chebyshev coefficients via the discrete
  cosine transform.

  \item \emph{Split.}
  Replace each local patch that fails the resolution
  criterion at degree $N_{\max}$ by two overlapping children,
  each initialized at degree $N_{\max}$. For a parent
  interval of width $h$, prescribe the overlap as a fraction
  $\eta$ of each child's width, giving
  \[
    h_{\mathrm{child}}=\frac{h}{2-\eta},
    \qquad
    h_{\mathrm{ov}}=\eta h_{\mathrm{child}}.
  \]
  The children cover the parent interval. Adjust neighboring
  overlaps as needed to preserve the requirement that only
  immediate neighbors overlap.

  \item \emph{Merge.}
  Sweep from left to right through the local patch chain.
  For each neighboring pair, construct a degree-$N_{\max}$
  interpolant of the current composite field~\eqref{eq:pou-blend} on their union.
  If it satisfies the resolution criterion for every field,
  replace the pair by a single patch. Repeat the sweep until
  no further merges are accepted or a prescribed maximum
  number of passes is reached.

  \item \emph{Chop.}
  Examine every resulting local patch, including those left
  unchanged, and reduce its polynomial degree according to
  the coefficient criterion. This step changes resolution
 without altering patch geometry.

  \item \emph{Repeat or terminate.}
  If the patch geometry has changed, return to the optimization
  step on the updated discretization. Otherwise, terminate
  the adaptive cycle, performing a final optimization if
  the chop pass has changed the polynomial degrees.
\end{enumerate}
For problems with multiple fields, we use the same patch boundaries
for all fields to simplify the residual grids used during
training. Thus, a patch is split
if any field is unresolved, and neighboring patches are merged
only if every field is resolved in their union. The chop pass, however, selects the polynomial degree independently
for each field.

This paper does not use adaptive refinement in the two-dimensional
experiments. We note that  the same coefficient-based approach
could be extended to higher dimensions. One option is to maintain
a tree of nonoverlapping rectangular cells and enlarge each
cell to form overlapping patches used throughout
training~\cite{aiton2019adaptive}. Alternatively, training
could use a nonoverlapping quadtree decomposition with
interface matching conditions~\cite{toselli2004domain},
followed by an extension of the local approximations to
overlapping patches for partition-of-unity blending to obtain global smoothness.
Nevertheless, for the 2D problem considered in this work we were able to  obtain highly accurate results without deploying this adaptive refinement.

\section{Numerical Results}\label{sec:results}
This section presents results for several self-similar blowup
problems: the inviscid Burgers equation,
the \ac{CCF} and \ac{CLM} models in the generalized De Gregorio family,
and the two-dimensional Boussinesq equations. The supplementary
materials also include preliminary results for the Euler equations.

For the problems treated with adaptivity, the local patch
decomposition and polynomial degrees are selected automatically.
For each field $q(y)$, we must still prescribe the following
parameters introduced in Section~\ref{sec:cheb-background}:
\begin{itemize}
  \item The leading-order far-field exponent $\sigma$ and
        associated factor $E_q(y)$ in~\eqref{eq:qhat}.
  \item The mapping exponent $\alpha$ of the asymptotic patch.
  \item The reference point $y_{\mathrm{mid}}^{\mathrm{asym}}$.
  \item The number of interpolation nodes $n_{\mathrm{asym}}$
        on the asymptotic patch.
  \item The extent of the local region, $y^{(m)}_\mathrm{end}$.
\end{itemize}
For training, we also prescribe two parameters introduced
in Section~\ref{sm-sec:loss}:
\begin{itemize}
  \item The residual-grid refinement factor $\kappa$.
  \item The gradient-penalty weight $c_{\mathrm{fgrad}}$.
\end{itemize}
We exploit the even or odd symmetry of each field about the
origin in the one-dimensional problems and about the $y_2$-axis
in the two-dimensional Boussinesq problem. Training is performed
on the corresponding half-interval or quarter plane, with symmetry
boundary conditions imposed at the symmetry interface. After
optimization, the fields are extended to the full domain
according to their prescribed parity and combined into a
smooth composite approximation, including across the symmetry
interface. The extension procedure is described
in Section~\ref{sm-sec:parity-extension}.

\begin{remark}[Grid selection]
  For brevity, we defer the derivation of the leading-order
  far-field behavior and mapping exponents, together with
  the grid-selection details and final discretizations, to
  Sections~\ref{sm-sec:exponents}--\ref{sm-sec:boussinesq} of the supplementary materials and to the code.
  For the two-dimensional Boussinesq problem, we use one local
  patch and one asymptotic patch without adaptivity and
  therefore specify both grids in full.
\end{remark}

\begin{remark}[Validation]\label{rem:val}
  All errors reported in this section are evaluated at $10^6$
  random points drawn according to the sampling procedure
  specified for each validation region. We report the sampled
  $\ell^\infty$ error and the root-mean-square error (RMSE).
  For the Burgers problem, whose profiles grow without bound,
  we instead report relative errors in certain regions, normalizing the maximum
  absolute error by the maximum absolute reference value and
  the RMSE by the root-mean-square reference value. The latter
  is equivalent to the relative $\ell^2$ error.
\end{remark}

\begin{remark}[Code and computer system]
  The code used to produce these results~\cite{blowuprepo} will be  publicly
available upon acceptance of the manuscript.
All reported computation times were measured on a single
NVIDIA A100 GPU with 80\,GB of memory.
\end{remark}

\subsection{The Burgers equation}\label{sec:burgers}
The inviscid Burgers equation
\begin{equation}\label{eq:burgers}
  u_t + u\,u_x = 0
\end{equation}
is  one of the simplest settings for finite-time gradient
blowup from smooth initial data. Its singularity formation
and self-similar structure are understood
analytically~\cite{eggers2008role}, making it a useful benchmark
for assessing the accuracy of our method. Taking the
self-similar ansatz with blowup time normalized to $t=1$,
\[
  u(x,t) = (1-t)^{\lambda}\,U(y),
  \qquad
  y = \frac{x}{(1-t)^{1+\lambda}},
  \qquad \lambda > 0,
\]
reduces~\eqref{eq:burgers} to the stationary profile equation
\begin{equation}\label{eq:ss-burgers}
  f[U](y)
  = -\lambda U
    + \bigl((1+\lambda)y + U\bigr)\,\partial_y U
  = 0,
  \qquad y\in\mathbb{R}.
\end{equation}
We seek odd profiles and therefore solve on the half-line
$y\geq 0$. The corresponding family of solutions is given
implicitly by
\begin{equation}\label{eq:burgers-exact}
  y = -U - C\operatorname{sgn}(U)\,|U|^{1+1/\lambda},
  \qquad C>0.
\end{equation}
These profiles are smooth across the origin precisely when
$1+1/\lambda$ is an odd integer, giving the discrete sequence
\begin{equation}\label{eq:burgers-smooth-lam}
  \lambda = \frac{1}{2i+2},
  \qquad i=0,1,2,\dots.
\end{equation}

We use the core loss components described in
Section~\ref{sm-sec:loss}, with $\kappa=3$ and
$c_{\mathrm{fgrad}}=1$, and enforce the following constraints
directly through fixed nodal values:
\begin{itemize}
  \item \emph{Symmetry boundary condition.}
  We impose $U(0)=0$, as required by odd symmetry.
  The first node of the innermost local patch lies at the
  origin; its value is fixed to zero and excluded from
  the optimization variables.

  \item \emph{Amplitude normalization.}
  We impose $\widehat U(\infty)=-1$, where $\widehat U$
  is the normalized field~\eqref{eq:qhat}.
  This fixes $C=1$ in~\eqref{eq:burgers-exact}.
  The endpoint node of the asymptotic patch corresponding
  to infinity is fixed to $-1$ and likewise excluded from
  the optimization variables.
\end{itemize}
No additional penalty terms are used for this problem.

Table~\ref{tab:burgers-results} reports results for
$\lambda=1/2$, $1/4$, and $1/6$, evaluated using the validation
procedure in Remark~\ref{rem:val}. We evaluate these quantities using $10^6$ samples drawn
uniformly in $y$ within the local region and a further
$10^6$ samples drawn via $y=\sinh(z)$, with $z$ uniform
on $[-40,40]$, for an ``infinite domain'' assessment (corresponding approximately to $y\in[-1.18\times10^{17},\,1.18\times10^{17}]$).
Both the equation residuals and the errors relative to the
exact profiles approach accuracies limited by floating-point precision.
For comparison, we implement a \ac{PINN}-based approach guided by the training framework of~\cite{wang2025discovery}. Table~\ref{tab:burgers-comparison} compares the computational cost, maximum residuals, and errors over the same ``infinite domain'' validation points. In all three cases, \ac{ODIL} achieves smaller errors relative to the exact solution while requiring substantially fewer parameters and several orders of magnitude less wall-clock time for training iterations than \acp{PINN}. Additional implementation details and results for several neural network architectures are provided in the supplementary materials (Supplementary Section~\ref{sm-sec:pinns}).
These results demonstrate that the method accurately and efficiently
recovers several smooth self-similar branches at prescribed
scaling exponents, providing a benchmark for the problems
below, where exact profiles are unavailable.
\begin{table}[htb]
  \centering
  \small
  \setlength{\tabcolsep}{5pt}
  \caption{Residuals and profile errors for the Burgers
    equation, evaluated at $10^6$ points per validation
    region. Local-region samples are uniform in $y$;
    ``infinite domain'' samples use $y=\sinh(z)$ with
    $z$ uniform on $[-40,40]$. Residuals are reported
    in absolute terms in both regions. Profile errors
    are absolute in the local region and relative
    in the ``infinite domain'', with RMSE and maximum
    error normalized by the root-mean-square and maximum
    absolute reference values, respectively.}
  \label{tab:burgers-results}
  \begin{tabular}{
    ll
    *{4}{S[table-format=1.2e-2]}
  }
    \toprule
    & & \multicolumn{2}{c}{$|y|\le 3$}
      & \multicolumn{2}{c}{``Infinite domain''} \\
    \cmidrule(lr){3-4}\cmidrule(lr){5-6}
    Profile & Quantity
      & {RMSE} & {Maximum}
      & {RMSE} & {Maximum} \\
    \midrule
    $\lambda=1/2$
      & $f$
      & 1.79e-14 & 4.47e-14
      & 3.69e-14 & 5.64e-14 \\
      & $U-U_{\mathrm{ex}}$
      & 1.06e-15 & 3.85e-15
      & 7.14e-16 & 1.23e-15 \\
    \addlinespace
    $\lambda=1/4$
      & $f$
      & 5.86e-15 & 2.27e-14
      & 5.73e-15 & 1.60e-14 \\
      & $U-U_{\mathrm{ex}}$
      & 5.71e-16 & 1.85e-15
      & 4.29e-16 & 8.96e-16 \\
    \addlinespace
    $\lambda=1/6$
      & $f$
      & 3.65e-15 & 1.97e-14
      & 2.71e-15 & 1.62e-14 \\
      & $U-U_{\mathrm{ex}}$
      & 3.20e-16 & 1.38e-15
      & 2.40e-16 & 7.59e-16 \\
    \bottomrule
  \end{tabular}
\end{table}

\begin{table}[t]
\centering
\caption{Comparison of \ac{ODIL} and \acp{PINN} for the Burgers
self-similar profiles. As in Table~\ref{tab:burgers-results}, the residuals $f$ are reported in absolute terms while the profile errors $U - U_\mathrm{ex}$ are relative.}
\label{tab:burgers-comparison}

\small
\renewcommand{\arraystretch}{1.3}
\setlength{\tabcolsep}{6pt}

\begin{tabular}{@{}c l r r r r@{}}
\toprule
$\lambda$
& Method
& Parameters (Stage 1/2)
& Time (s)
& Max. $f$
& Max. $U - U_\mathrm{ex}$ \\
\midrule

$1/2$
& ODIL
& 84/--
& 0.405
& $5.64\times10^{-14}$
& $1.23\times10^{-15}$ \\

& PINN
& 321/817
& 159
& $9.22\times10^{-15}$
& $4.52\times10^{-15}$ \\

\midrule

$1/4$
& ODIL
& 109/--
& 0.620
& $1.60\times10^{-14}$
& $8.96\times10^{-16}$ \\

& PINN
& 673/1,801
& 504
& $5.46\times10^{-14}$
& $2.81\times10^{-15}$ \\

\midrule

$1/6$
& ODIL
& 127/--
& 0.754
& $1.62\times10^{-14}$
& $7.59\times10^{-16}$ \\

& PINN
& 673/1,801
& 872
& $8.93\times10^{-14}$
& $5.40\times10^{-15}$ \\

\bottomrule
\end{tabular}%

\end{table}
\begin{remark}[Inference of $\lambda$]\label{rem:burgers-lambda}
  Although this paper focuses on computing profiles at prescribed
  values of $\lambda$, the method also offers a possible approach
  to inferring the admissible $\lambda$ values similarly  to  previous works~\cite{wang2023asymptotic,wang2025discovery}. 
  For each candidate
  $\lambda$, one can solve the profile problem, apply the
  parity-extension procedure described in Section~\ref{sm-sec:parity-extension},
  and use the residual near the origin to define the objective
  \begin{equation}\label{eq:burgers-symmetry-loss}
    \mathcal{L}(\lambda)
    = \frac{1}{N}\sum_{k=1}^{N}
      f[\widetilde U_\lambda](y_k)^2,
  \end{equation}
  where $\widetilde U_\lambda$ is the extended profile and the
  $y_k$ are sampled uniformly at random from
  $|y|\leq\delta y_{\mathrm{end}}^{(1)}$ for some prescribed positive parameter $\delta$ close to 0.

 Here, we minimize this objective
  using bounded Brent minimization~\cite{brent2013algorithms}
  over $\lambda\in[1/(2i+3),1/(2i+1)]$ to isolate the $i$th branch,
  rebuilding the grids and warm-starting each profile solve
  from the previous solution. With $N=2000$ and a search
  tolerance of $10^{-12}$, we recover $\lambda=1/2$, $1/4$,
  and $1/6$ with absolute errors of $1.8\times10^{-9}$,
  $1.1\times10^{-8}$, and $4.7\times10^{-7}$, respectively.
  The maximum sampled residuals of the extended profiles
  are of order $10^{-13}$ or below. A systematic study of
  $\lambda$ inference is a subject of future work.
\end{remark}

\subsection{The generalized De Gregorio family}\label{sec:degregorio}
We consider the one-parameter family of one-dimensional models
introduced in~\cite{okamoto2008generalization},
\begin{equation}\label{eq:dg-family}
  \omega_t + a\,u\,\omega_x = \omega\,u_x,
  \qquad
  u_x = \Hop[\omega],
\end{equation}
posed on the real line, where $\Hop$ denotes the Hilbert transform,
\begin{equation}\label{eq:hilbert-def}
  \Hop[\omega](x)
  = \frac{1}{\pi}\,\mathrm{p.v.}\!\int_{-\infty}^{\infty}
    \frac{\omega(x')}{x-x'}\,dx'.
\end{equation}
We study two members of this family: the \ac{CLM}
model~\cite{constantin1985simple}, corresponding to $a=0$,
for which advection is absent and
$\omega_t=\omega\,\Hop[\omega]$; and the \ac{CCF}
model~\cite{cordoba2005formation}, corresponding to $a=-1$.
The family provides a tractable setting for studying how
nonlocal stretching and advection influence singularity
formation, with $a$ controlling the strength and direction
of advection~\cite{elgindi2021finite}. Because $a$ appears
in only one term, the same numerical formulation applies
throughout the family. Here, we report results for the
\ac{CLM} model and three self-similar branches of the
\ac{CCF} model.

Taking the backward self-similar ansatz with blowup time
normalized to $t=1$,
\begin{equation}\label{eq:dg-ansatz}
  \omega(x,t) = \frac{\Omega(y)}{1-t},
  \qquad
  u(x,t) = (1-t)^{\lambda}\,U(y),
  \qquad
  y = \frac{x}{(1-t)^{1+\lambda}},
\end{equation}
reduces~\eqref{eq:dg-family} to the stationary system
\begin{align}
  f_1[\Omega,U](y)
  &= \Omega
    + \bigl((1+\lambda)y+a\,U\bigr)\,\partial_y\Omega
    - \Omega\,\partial_y U = 0,
    \label{eq:dg-f1}\\
  f_2[\Omega,U](y)
  &= \partial_y U-\Hop[\Omega] = 0.
    \label{eq:dg-f2}
\end{align}
We seek odd profiles $\Omega$. Since the Hilbert transform
of an odd function is even, equation~\eqref{eq:dg-f2}
implies that $\partial_y U$ is even. Choosing $U(0)=0$
then makes $U$ odd as well.

We use the core loss components described in
Section~\ref{sm-sec:loss}, with $\kappa=3$ and
a problem-dependent gradient-penalty weight
$c_{\mathrm{fgrad}}$, together with the following constraints:
\begin{itemize}
  \item \emph{Symmetry boundary conditions.}
  We impose $\Omega(0)=U(0)=0$ by fixing both field values
  at the first node of the innermost local patch and
  excluding them from the optimization variables.

  \item \emph{Normalization.}
  We rescale the profiles using either
  a penalty on $\Omega$ or $\partial_y\Omega$ weighted by a prescribed factor 
  $\sqrt{c_g}$, following~\cite{wang2023asymptotic}, or a prescribed
  value of $\widehat\Omega$ at infinity, where
  $\widehat\Omega$ is the normalized field~\eqref{eq:qhat}.
  In the latter case, the corresponding endpoint value
  on the asymptotic patch is fixed and excluded from
  the optimization variables. The specific normalization
  is given for each problem below.
\end{itemize}
We use a two-stage optimization to improve convergence.
In the first stage, normalization is enforced through a
soft penalty and $c_{\mathrm{fgrad}}>0$. In the second stage,
we fix $\widehat\Omega(\infty)$ at its learned value from the first stage and
set $c_{\mathrm{fgrad}}=0$. This second stage was
empirically helpful in reducing the equation residuals
to levels near machine precision.

For the \ac{CLM} model ($a=0$) with $\lambda=0$, the
stationary system admits the exact profiles
\begin{equation}\label{eq:clm-exact}
  \Omega_{\mathrm{ex}}(y) = -\frac{2y}{1+y^2},
  \qquad
  U_{\mathrm{ex}}(y) = 2\arctan(y).
\end{equation}
Both profiles are odd, and
$\partial_y\Omega_{\mathrm{ex}}(0)=-2$.
Moreover,
\begin{equation}\label{eq:clm-hilbert-exact}
  \Hop[\Omega_{\mathrm{ex}}](y)
  = \partial_y U_{\mathrm{ex}}(y)
  = \frac{2}{1+y^2}.
\end{equation}
These identities provide accuracy tests for both the computed
profiles and the numerical evaluation of the Hilbert transform.
For the \ac{CCF} model, we optimize the profiles at three prescribed
values of $\lambda$ obtained in previous
work~\cite{wang2025discovery,wang2025resolving}: $\lambda=1.180777662899$
for the stable branch, and $\lambda=0.6057337012032$ and
$\lambda=0.471324863862$ for the first two unstable branches.

\begin{figure}[htb]
  \centering
  \includegraphics[width=\textwidth]{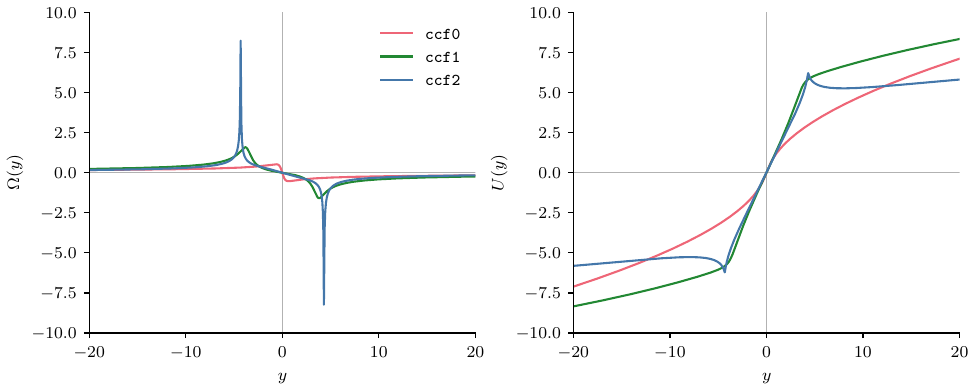}
  \caption{Computed \ac{CCF} profiles on $[-20,20]$, extended to $y<0$ by
    the parity-extension procedure of Section~\ref{sm-sec:parity-extension}:
    $\Omega$ (left) and $U$ (right) for \texttt{ccf0},
    \texttt{ccf1} and \texttt{ccf2}. The narrow interior feature of
    \texttt{ccf2} near $|y| \approx 4.2$ drives both its patch count and
    its hand-placed quadrature segments.}
  \label{fig:dg-profiles}
\end{figure}

\begin{figure}[htb]
  \centering
  \includegraphics[width=\textwidth]{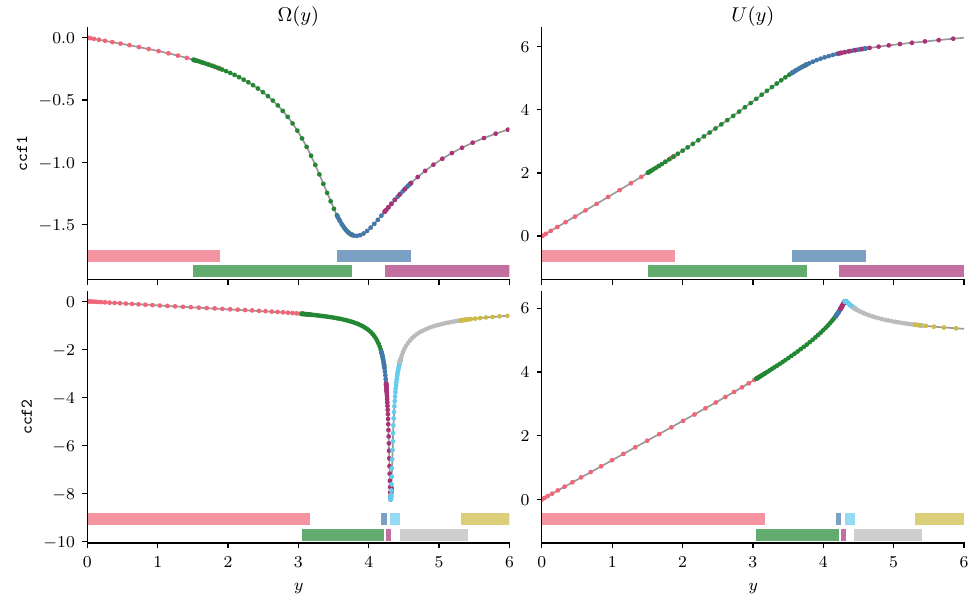}
  \caption{Local-patch decomposition of \texttt{ccf1} (top) and
    \texttt{ccf2} (bottom) on $[0,6]$, for $\Omega$ (left) and $U$
    (right). Markers are the interpolation nodes of each patch, colored by
    patch; the bars below each curve give the patch extents, staggered to
    show their overlaps. Note that the adaptive procedure can find features of the profiles.}
  \label{fig:dg-patches}
\end{figure}

Table~\ref{tab:clm-ccf-results} reports the equation residuals
and, for \ac{CLM}, the errors between the computed and exact profiles.
To make residual computations fair, points extremely close to the Hilbert-transform quadrature nodes are
excluded to avoid numerical cancellation that is artificially inflating residual metrics. The residuals for \ac{CLM} and the stable and first unstable
\ac{CCF} branches reach levels limited by machine epsilon,
with the exact \ac{CLM} solution providing an additional
check on accuracy: the maximum error in $\Omega$ is
approximately $1.76\times10^{-12}$, within roughly one
order of magnitude of the maximum residuals. The second
unstable \ac{CCF} branch also achieves high accuracy,
with maximum residuals of order $10^{-10}$, although
it does not reach the same precision as the other branches;
a possible explanation is discussed in
Remark~\ref{rem:ill-conditioning}.
Figure~\ref{fig:dg-profiles} shows the computed \ac{CCF} profiles, and
Figure~\ref{fig:dg-patches} shows the local-patch decompositions of the
two unstable branches in the near field, where the narrow interior
feature of the second unstable branch near $|y|\approx 4.2$ is resolved
by three narrow patches. This accuracy is obtained at a modest computational cost.
For the first unstable \ac{CCF} branch, optimization at the
prescribed $\lambda$ required 48 stage-1 \ac{LM} iterations and 50 stage-2 \ac{LM} iterations (13 seconds total),
whereas~\cite{wang2025discovery} reports three hours on an
A100 GPU to reach residuals of order $10^{-8}$. Reducing computation times from hours to tens of seconds allows
ideas to be tested and refined much more quickly, supporting
both numerical method development and the search for
self-similar solutions in more challenging problems.

\begin{table}[htb]
  \centering
  \small
  \setlength{\tabcolsep}{5pt}
  \caption{Residuals and profile errors for the \ac{CCF} and
  \ac{CLM} models, evaluated at $10^6$ points per validation
  region. Local-region samples are uniform in the original
  coordinate $y$. The ``infinite domain'' columns use
  $y=\sinh(z)$ with $z$ sampled uniformly on $[-40,40]$.
  Points close to Hilbert-transform quadrature nodes are
  excluded to avoid numerical cancellation in the quadrature
  evaluation. Columns report the RMSE and maximum absolute
  value of each residual or profile error over the samples.}
  \label{tab:clm-ccf-results}
  \begin{tabular}{
    ll
    *{4}{S[table-format=1.2e-2]}
  }
    \toprule
    & & \multicolumn{2}{c}{$|y|\le y^{(m)}_\mathrm{end}$}
      & \multicolumn{2}{c}{``Infinite domain''} \\
    \cmidrule(lr){3-4}\cmidrule(lr){5-6}
    Model & Quantity
      & {RMSE} & {Maximum}
      & {RMSE} & {Maximum} \\
    \midrule
    \texttt{clm}
    & $f_1$
    & 1.23e-14 & 7.69e-14
    & 3.66e-14 & 1.70e-13 \\
    & $f_2$
    & 1.85e-14 & 1.23e-13
    & 7.16e-15 & 1.01e-13 \\
    & $\Omega-\Omega_{\mathrm{ex}}$
    & 9.07e-13 & 1.76e-12
    & 2.71e-13 & 1.76e-12 \\
    \midrule
    \texttt{ccf0} ($\lambda=1.180777662899$)
      & $f_1$
      & 3.29e-15 & 4.54e-14
      & 2.00e-15 & 3.38e-14 \\
      & $f_2$
      & 5.57e-15 & 9.09e-14
      & 6.18e-15 & 5.43e-14 \\
    \addlinespace
    \texttt{ccf1} ($\lambda=0.6057337012032$)
      & $f_1$
      & 2.15e-14 & 4.03e-13
      & 2.10e-14 & 2.71e-13 \\
      & $f_2$
      & 4.37e-14 & 3.05e-13
      & 1.92e-13 & 7.22e-13 \\
    \addlinespace
    \texttt{ccf2} ($\lambda=0.471324863862$)
      & $f_1$
      & 3.78e-12 & 1.04e-10
      & 7.34e-12 & 1.04e-10 \\
      & $f_2$
      & 3.74e-12 & 1.37e-10
      & 5.55e-12 & 1.36e-10 \\
    \bottomrule
  \end{tabular}
\end{table}

\begin{remark}[Ill-conditioning]\label{rem:ill-conditioning}
  Resolving the sharp features of the second unstable
  \ac{CCF} profile requires high polynomial degrees and
  narrow overlaps between patches, both of which can
  amplify errors during differentiation. On a patch of
  width $h$, the first-derivative operator for degree-$n$
  Chebyshev polynomials has norm $O(n^2/h)$, while
  differentiating the partition-of-unity weights introduces
  factors proportional to the inverse overlap
  width~\cite{trefethen2013approximation,aiton2017adaptive}.
  These effects may contribute to ill-conditioning and
  limit the attainable residual accuracy. One possible remedy is to represent field derivatives
  as the primary unknowns and recover the fields by
  integration. Such an integral
  reformulation can substantially improve the conditioning
  of spectral discretizations and would be straightforward to implement in the \ac{ODIL} framework.
  Adapting this approach to our formulation could permit
  higher polynomial degrees and wider overlaps while
  reducing the amplification associated with direct
  differentiation.
\end{remark}

\FloatBarrier
\subsection{The two-dimensional Boussinesq equations}\label{sec:boussinesq}
On the half-plane $x_2 \ge 0$, the inviscid Boussinesq system is
\begin{equation}\label{eq:boussinesq}
  \partial_t u + u\cdot\nabla u + \nabla p = (0,\theta),
  \qquad
  \partial_t \theta + u\cdot\nabla \theta = 0,
  \qquad
  \nabla\cdot u = 0,
\end{equation}
where $u = (u_1,u_2)$ is the velocity and $\theta$ is the scalar
temperature, subject to the nonpenetration condition
\begin{equation}\label{eq:nonpenetration}
  u_2(x_1,0,t) = 0 .
\end{equation}
Under axial symmetry and up to exponentially small corrections, this
system is closely connected to the Luo--Hou blowup scenario for the
three-dimensional Euler equations in the presence of a cylindrical
boundary~\cite{luo2014potentially,luo2014toward,wang2023asymptotic}. The
self-similar profile associated with that scenario was first identified
numerically using \acp{PINN}~\cite{wang2023asymptotic} and later
established rigorously by computer-assisted
proof~\cite{chen2022stable,chen2025stable}.

Following~\cite{wang2023asymptotic}, we take the backward self-similar
ansatz with blowup time normalized to $t=1$,
\begin{equation}\label{eq:selfsimilar-ansatz}
  u(x,t) = (1-t)^{\lambda}\,U(y),
  \qquad
  \theta(x,t) = (1-t)^{-1+\lambda}\,\Theta(y),
  \qquad
  y = \frac{x}{(1-t)^{1+\lambda}},
\end{equation}
with scaling exponent $\lambda > -1$, under which~\eqref{eq:boussinesq}
becomes the stationary system
\begin{equation}\label{eq:ss-boussinesq}
  \begin{aligned}
    -\lambda U + \bigl((1+\lambda)y + U\bigr)\cdot \nabla U + \nabla P &= (0,\Theta),\\
    (1-\lambda)\Theta + \bigl((1+\lambda)y + U\bigr)\cdot \nabla \Theta &= 0,\\
    \nabla\cdot U &= 0 .
  \end{aligned}
\end{equation}
Solutions of~\eqref{eq:ss-boussinesq} are expected to have infinite
energy~\cite{chae2007nonexistence}; only mild growth at infinity is
required, which permits such a solution to be cut off to produce a
finite-energy asymptotically self-similar solution.

\begin{figure}[htb]
  \centering
  \includegraphics[width=0.85\textwidth]{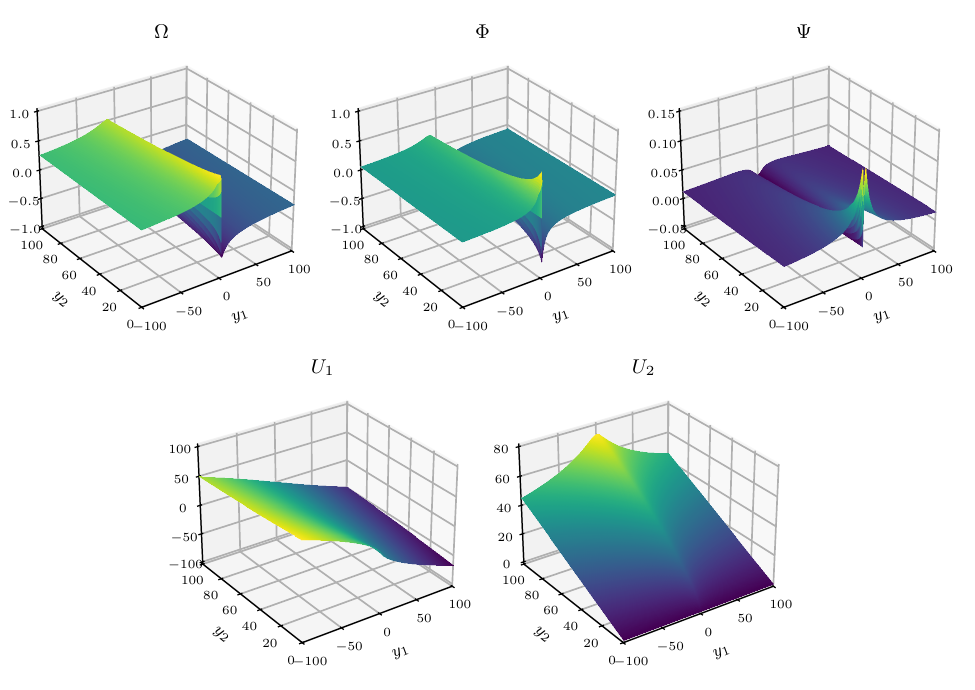}
  \caption{The five computed fields of the two-dimensional Boussinesq
    profile on the near field $[-100,100]\times[0,100]$, extended to
    $y_1<0$ with the parities~\eqref{eq:symmetry} by the procedure of
    Section~\ref{sm-sec:parity-extension}.}
  \label{fig:bq-fields}
\end{figure}

Introducing the vorticity and the temperature gradient,
\begin{equation}\label{eq:omega-phi-psi}
  \Omega = \partial_{y_1}U_2 - \partial_{y_2}U_1,
  \qquad
  \Phi = \partial_{y_1}\Theta,
  \qquad
  \Psi = \partial_{y_2}\Theta,
\end{equation}
and the self-similar transport field
\begin{equation}\label{eq:transport-field}
  V = (1+\lambda)\,y + U,
\end{equation}
we solve~\eqref{eq:ss-boussinesq} in vorticity form for the five unknown
fields
\begin{equation}\label{eq:bq-unknowns}
  \Omega, \quad \Phi, \quad \Psi, \quad U_1, \quad U_2 ,
\end{equation}
by driving the six residuals
\begin{equation}\label{eq:bq-residuals}
  \begin{aligned}
    f_1 &= \Omega + V\cdot\nabla\Omega - \Phi,
    & f_4 &= \partial_{y_1}U_1 + \partial_{y_2}U_2, \\
    f_2 &= \bigl(2+\partial_{y_1}U_1\bigr)\Phi + V\cdot\nabla\Phi
            + (\partial_{y_1}U_2)\,\Psi,
    & f_5 &= \Omega + \bigl(\partial_{y_2}U_1 - \partial_{y_1}U_2\bigr), \\
    f_3 &= \bigl(2+\partial_{y_2}U_2\bigr)\Psi + V\cdot\nabla\Psi
            + (\partial_{y_2}U_1)\,\Phi,
    & f_6 &= \partial_{y_2}\Phi - \partial_{y_1}\Psi
  \end{aligned}
\end{equation}
to zero. The residuals $f_1$, $f_2$ and $f_3$ are the transport
equations for $\Omega$, $\Phi$ and $\Psi$, $f_4$ is the divergence-free
condition, $f_5$ is the definition of the vorticity, and $f_6$ is the
integrability condition that makes $(\Phi,\Psi)$ a gradient field. The
velocity components are treated as primitive unknowns, so the
velocity--vorticity relation is enforced through $f_4$ and $f_5$ rather
than by inverting an elliptic operator; further discussion of the
formulation is given in Section~\ref{sm-sec:bq-formulation}.

\begin{figure}[htb]
  \centering
  \includegraphics[width=0.85\textwidth]{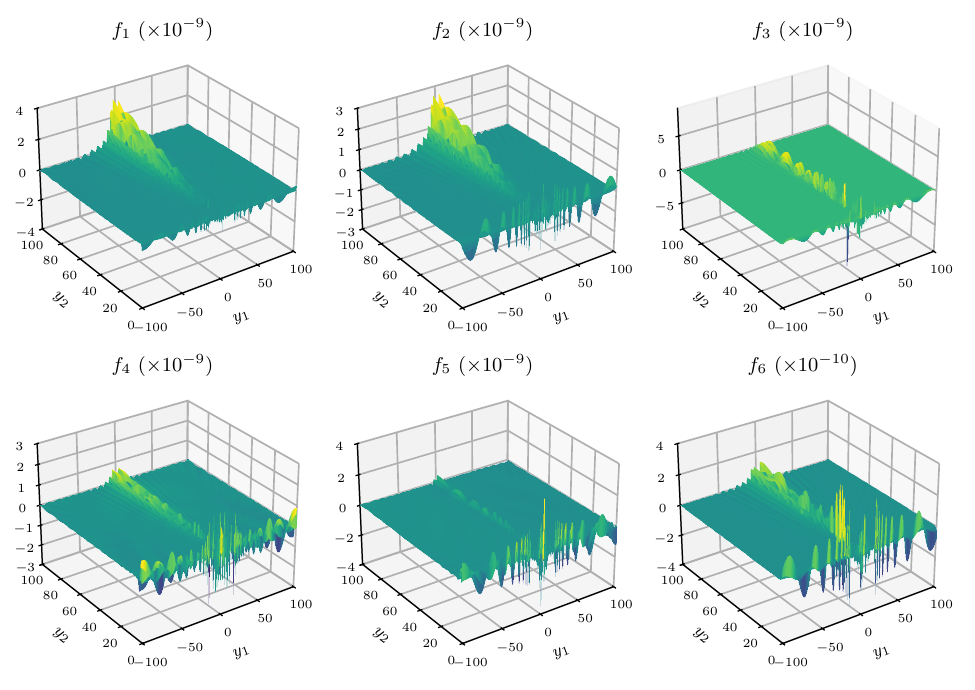}
  \caption{The six residuals~\eqref{eq:bq-residuals} of the
    two-dimensional Boussinesq profile on
    $[-100,100]\times[0,100]$, evaluated on the parity-extended fields of
    Figure~\ref{fig:bq-fields}. Each residual is scaled by the power of ten
    given in its panel title.}
  \label{fig:bq-residuals}
\end{figure}

As in~\cite{wang2023asymptotic}, we impose the symmetry conditions
\begin{equation}\label{eq:symmetry}
  (U_1,\Phi,\Omega)\ \text{odd in } y_1,
  \qquad
  (U_2,\Psi)\ \text{even in } y_1,
\end{equation}
together with the wall condition inherited from~\eqref{eq:nonpenetration},
\begin{equation}\label{eq:ss-boundary}
  U_2(y_1,0) = 0 ,
\end{equation}
and the normalization
\begin{equation}\label{eq:normalization}
  \partial_{y_1}\Omega(0,0) = -1 ,
\end{equation}
which removes the scaling degeneracy of the self-similar system. By
symmetry, the computation is carried out on the quarter plane
\begin{equation}\label{eq:quarter-plane}
  \{\,y = (y_1,y_2) : y_1 \ge 0,\ y_2 \ge 0\,\}.
\end{equation}

We use the core loss components described in
Section~\ref{sm-sec:loss}, with $\kappa=2$ and
$c_{\mathrm{fgrad}}=1$, together with the following constraints:
\begin{itemize}
  \item \emph{Odd symmetry.}
  For $\Omega$, $\Phi$ and $U_1$, we impose $q(0,y_2)=0$ by fixing the
  nodal values on the symmetry axis of both patches to zero and
  excluding them from the optimization variables.

  \item \emph{Even symmetry.}
  For $\Psi$ and $U_2$, we impose the homogeneous Neumann condition
  $\partial_{y_1}q(0,y_2)=0$ as a penalty evaluated on the symmetry axis
  of both patches, weighted by $\sqrt{c_{\mathrm{bc}}}$.

  \item \emph{Wall condition.}
  We impose~\eqref{eq:ss-boundary} by fixing the nodal values of $U_2$
  on the wall to zero and excluding them from the optimization
  variables.

  \item \emph{Normalization.}
  We impose~\eqref{eq:normalization} through the scalar penalty
  $\partial_{y_1}\Omega(0,0)+1$, weighted by $\sqrt{c_g}$.
\end{itemize}
We use $c_{\mathrm{bc}}=c_g=1$.

Each field is represented on one Cartesian local patch covering
$[0,15]^2$ and one polar asymptotic patch starting at
$r_{\mathrm{start}}=10$, without adaptivity (see Section~\ref{sec:methodology-general-2d}). The far-field exponents and
the mapping exponents are derived in
Section~\ref{sm-sec:exponents}, and the discretizations are given in
Section~\ref{sm-sec:bq-optimize} and Table~\ref{sm-tab:bq-cfg}. 

Table~\ref{tab:bq-results} reports the equation residuals for the prescribed
scaling exponent $\lambda = 1.9205599746927$ obtained
in~\cite{wang2025discovery}, evaluated
using $10^6$ samples drawn uniformly by area in the half-disk
$\{|y|\le 15,\ y_2\ge 0\}$ and a further $10^6$ samples drawn via
$y_i=\sinh(z_i)$, with $z_1$ uniform on $[-20,20]$ and $z_2$ uniform on
$[0,20]$, for an ``infinite domain'' assessment extending to
$|y|\approx 3.4\times10^{8}$. The fields are extended to $y_1<0$ by the
parity-extension procedure of Section~\ref{sm-sec:parity-extension}. In
the near field, all six residuals have maximum values below
$1.2\times10^{-8}$ and RMSE below $4\times10^{-10}$. Over the infinite
domain, the maximum residuals are at most $1.0\times10^{-7}$, attained
by $f_5$, so the accuracy is largely retained in the far field although
no condition is imposed there.
Figures~\ref{fig:bq-fields} and~\ref{fig:bq-residuals} show the computed
fields and the residuals in the near field. Within the
local region, our method achieves residual levels comparable
to those of the \ac{PINN}-based approach
in~\cite{wang2025discovery}, with limited tuning of the
discretization and no adaptive refinement. These promising results motivate extending adaptive refinement
to two dimensions, with the goal of computing high-precision
self-similar profiles for open problems such as singularity
formation in the unforced three-dimensional Euler
equations~\cite{ganeshram2026stable}. Preliminary results regarding this problem
are presented in Supplementary Section~\ref{sm-sec:euler} using the Chebyshev \ac{ODIL} framework but on a truncated domain.

\begin{table}[htb]
  \centering
  \small
  \setlength{\tabcolsep}{5pt}
  \caption{Residuals for the two-dimensional Boussinesq equations,
    evaluated at $10^6$ points per validation region. Near-field samples
    are uniform by area in the half-disk $\{|y|\le 15,\ y_2\ge 0\}$. The
    ``infinite domain'' columns use $y_i=\sinh(z_i)$ with $z_1$ sampled
    uniformly on $[-20,20]$ and $z_2$ uniformly on $[0,20]$. Columns
    report the RMSE and maximum absolute value of each residual over the
    samples.}
  \label{tab:bq-results}
  \begin{tabular}{
    l
    *{4}{S[table-format=1.2e-2]}
  }
    \toprule
    & \multicolumn{2}{c}{$|y|\le 15$}
      & \multicolumn{2}{c}{``Infinite domain''} \\
    \cmidrule(lr){2-3}\cmidrule(lr){4-5}
    Residual
      & {RMSE} & {Maximum}
      & {RMSE} & {Maximum} \\
    \midrule
    $f_1$ & 9.27e-11 & 4.26e-09 & 6.88e-09 & 4.48e-08 \\
    $f_2$ & 5.62e-11 & 1.66e-09 & 2.32e-09 & 1.45e-08 \\
    $f_3$ & 3.70e-10 & 1.13e-08 & 3.95e-10 & 1.06e-08 \\
    $f_4$ & 1.91e-10 & 2.69e-09 & 2.24e-09 & 2.63e-08 \\
    $f_5$ & 2.88e-10 & 4.06e-09 & 1.63e-08 & 1.00e-07 \\
    $f_6$ & 3.69e-11 & 9.75e-10 & 2.25e-11 & 9.79e-10 \\
    \bottomrule
  \end{tabular}
\end{table}

\FloatBarrier
\section{Conclusions}\label{sec:conclusions}
We have presented computations of  self-similar singularity profiles on unbounded domains using an adaptive spectral formulation of \ac{ODIL}. Each
field is represented by its nodal values on overlapping local Chebyshev
patches joined by a smooth partition of unity, together with an
asymptotic patch whose grid is mapped from $[y_{\mathrm{start}},\infty]$
onto the reference interval. The leading far-field
power laws are factored out, while their amplitudes, including their
angular dependence in two dimensions, are determined by the optimization,
so that neither domain truncation nor prescribed far-field conditions are
required. The mapping exponent is selected from the asymptotic
expansion of the profile equations, which preserves the smoothness of the
normalized fields at infinity. Adaptive refinement guided by
Chebyshev-coefficient decay and damped Gauss--Newton iterations addresses
the representation and optimization errors within a common framework. The method computes the Burgers profiles to near machine precision, the
generalized De Gregorio profiles with residuals between $O(10^{-14})$ and
$O(10^{-10})$, and the two-dimensional Boussinesq profile with residuals
of $O(10^{-8})$; where timings are available for comparison, the cost is
substantially below that of \ac{PINN}-based approaches. We make the source code, solver configurations, computed profiles, and postprocessing scripts openly available to support independent reproduction and further development.
 
Several directions remain open. The attainable accuracy for profiles with
sharp interior features appears to be limited by the conditioning of
high-degree differentiation across narrow overlaps, which an integral
reformulation with field derivatives as primary unknowns may alleviate
(Remark~\ref{rem:ill-conditioning}). Adaptive refinement has so far been
applied only in one dimension, and its extension to two dimensions is the
natural next step towards high-precision profiles for open problems such
as singularity formation in the unforced three-dimensional Euler
equations~\cite{ganeshram2026stable}, for which preliminary results are
reported in Supplementary Section~\ref{sm-sec:euler}. While this paper focuses on fluid-related \acp{PDE}, this construction is motivated by a class of profile problems rather than a particular fluid equation. Backward self-similar heat profiles and thin-film rupture branches provide natural candidates for extension because they combine regular cores with algebraic far-field behavior~\cite{collot2019stability,dallaston2017self}. Dispersive profiles offer another direction, although their oscillatory asymptotics can require factoring a phase as well as an amplitude envelope~\cite{dahne2024self}. Such extensions require problem-specific asymptotic analysis; the present formulation targets stationary similarity profiles whose normalized fields are smooth in suitable mapped coordinates. More generally, the techniques presented in this paper, namely the nodal Chebyshev parameterization within \ac{ODIL}, overlapping patches with adaptive refinement, and an asymptotic patch whose far-field amplitudes are unknowns of the optimization, are not specific to self-similar blowup, and can be applied to a wide range of nonlinear problems for \acp{PDE}, in particular those posed on unbounded domains.

\bibliographystyle{unsrt}
\bibliography{refs}

\end{document}

\maketitle

\section{Methodology}\label{sm-sec:methodology}

Section~\ref{sec:methodology-general} of the main text describes the
multi-patch Chebyshev representation and the residual-based loss at the
level of detail needed to follow the results. This section supplies the
definitions omitted there. Section~\ref{sm-sec:cheb} collects the properties
of Chebyshev interpolation on which the method relies.
Section~\ref{sm-sec:grid} gives the map of the asymptotic patch and the
partition of unity that enter the one-dimensional composite
approximation~\eqref{eq:pou-blend}, and Section~\ref{sm-sec:2d} gives the
corresponding construction on the quarter plane. Section~\ref{sm-sec:loss}
specifies the blocks of the residual vector and the linear solvers used
within the Gauss--Newton iterations, and Section~\ref{sm-sec:adaptivity}
supplements the adaptive procedure of Section~\ref{sec:general-adaptive}.

\subsection{Chebyshev interpolation}\label{sm-sec:cheb}

The main text approximates each field on each patch by a Chebyshev
interpolant and uses the decay of its Chebyshev coefficients both to assess
resolution and to guide adaptive refinement
(Sections~\ref{sec:cheb-background} and~\ref{sec:general-adaptive}). We
recall here the definitions and results that underlie these uses. For
$n\ge2$, the interpolation nodes on $[-1,1]$ are the $n$
Chebyshev--Lobatto points
\begin{equation}\label{sm-eq:cheb-nodes}
  \xi^{(n)} = \left\{\xi_i^{(n)} =\cos\left(\frac{i\pi}{n-1}\right): 0\leq i \leq n-1\right\},
\end{equation}
and $q_n$ denotes the unique polynomial of degree at most $n-1$ that agrees
with $q$ at these points. The super-algebraic convergence invoked in the
main text is made precise by the following classical estimate, which shows
that the interpolation error decreases faster than any fixed power of $n$
when $q$ is smooth.

\begin{theorem}[{\cite[Thm.~7.2]{trefethen2013approximation}}]\label{sm-thm:cheb-convergence}
  Suppose $q^{(\nu)}$ has bounded variation $V_\nu$ on $[-1,1]$, with
  $q, q', \dots, q^{(\nu-1)}$ absolutely continuous there. Then for
  $n > \nu + 1$,
  \[
    \norm{q - q_n}_\infty \leq \frac{4V_\nu}{\pi \nu (n-1-\nu)^{\nu}}.
  \]
\end{theorem}

The coefficient form~\eqref{eq:cheb-expansion} of the interpolant
approximates the Chebyshev series of $q$,
\begin{equation}\label{sm-eq:cheb-coeffs}
  q(y) = \sum_{k=0}^\infty a_k T_k(y),
  \qquad
  a_k = \frac{2}{\pi}\int_{-1}^1 \frac{q(y)T_k(y)}{\sqrt{1-y^2}}\,dy
  \quad (k\ge1),
\end{equation}
where the coefficient $a_0$ carries an additional factor $\tfrac12$. The
decay of $|a_k|$ with $k$ reflects the smoothness of $q$, and
Section~\ref{sm-sec:adaptivity} makes this relation quantitative for the
resolution criterion of the adaptive procedure.

Derivatives are likewise computed in coefficient space. Writing
$q_n'(y) = \sum_{k=0}^{n-2} a_k^{(1)} T_k(y)$ and setting
$a_{n-1}^{(1)} = a_n^{(1)} = 0$, the coefficients of the derivative follow
from the downward recurrence
\begin{equation}\label{sm-eq:cheb-diff-recurrence}
  c_k\, a_k^{(1)} = a_{k+2}^{(1)} + 2(k+1)a_{k+1},
  \qquad k = n-2, n-3, \dots, 0,
\end{equation}
with $c_0=2$ and $c_k=1$ for $k \ge 1$~\cite{boyd2001chebyshev,trefethen2000spectral}.
Combining~\eqref{sm-eq:cheb-diff-recurrence} with the forward and inverse
coefficient transforms yields a differentiation operator on nodal values
that is exact for polynomials of degree at most $n-1$. All derivatives in
the residuals are evaluated in this way, and the chain and product rules
then account for the coordinate maps, the envelopes $E_q$ and the
partition-of-unity weights.

\subsection{One-dimensional multi-patch representation}\label{sm-sec:grid}

Section~\ref{sec:cheb-background} of the main text introduces the local
patches $K_i$, the asymptotic patch $K_{\mathrm{asym}}$ and the normalized
field~\eqref{eq:qhat}, and states that the asymptotic patch is mapped to the
reference interval by an algebraic map whose exponent $\alpha$ is chosen to
keep the normalized fields smooth. We now give this map explicitly and
describe the partition of unity that blends the patchwise interpolants
in~\eqref{eq:pou-blend}. Throughout, each patch overlaps only its immediate
neighbors.

Given the mapping exponent $\alpha > 0$ and a reference point
$y_{\mathrm{mid}}^\mathrm{asym} > y_{\mathrm{start}}^\mathrm{asym}$, we set
$L = \left(y_{\mathrm{mid}}^{\mathrm{asym}}\right)^{\alpha} -
\left(y_{\mathrm{start}}^{\mathrm{asym}}\right)^{\alpha}$ and define the
map~\cite{boyd2001chebyshev,boyd1987spectral,boyd1987orthogonal}
\begin{equation}\label{sm-eq:map}
  x(y) = \frac{\bigl(y^{\alpha} - \left(y_{\mathrm{start}}^{\mathrm{asym}}\right)^{\alpha}\bigr) - L}
              {\bigl(y^{\alpha} - \left(y_{\mathrm{start}}^{\mathrm{asym}}\right)^{\alpha}\bigr) + L},
\end{equation}
which sends $y = y_{\mathrm{start}}^\mathrm{asym}$ to $x=-1$,
$y = y_{\mathrm{mid}}^\mathrm{asym}$ to $x=0$, and $y \to \infty$ to $x=1$.
On the asymptotic patch, the normalized field is approximated by the
interpolant $\widehat q_{n_{\mathrm{asym}}}$ on an $n_{\mathrm{asym}}$-point
Chebyshev--Lobatto grid in $x$, evaluated at $x = x(y)$. Because a
polynomial is bounded on the closed reference interval, only a bounded
function can be approximated in this way; this is why the envelope $E_q$ is
factored out before interpolation whenever $q$ grows at infinity.

The exponent $\alpha$ controls the behavior of the map near the point at
infinity. Inverting the map~\eqref{sm-eq:map} gives
\[
  y(x)^\alpha = \left(y_{\mathrm{start}}^{\mathrm{asym}}\right)^\alpha + L\,\frac{1+x}{1-x}
\]
and thus
\[
  y(x) \sim \left(\frac{2L}{1-x}\right)^{1/\alpha} \quad \text{as } x \rightarrow 1^-.
\]
Consequently, for a normalized field $\widehat{q}(y)$ with the expansion
$\widehat{q}(y) = c_0 + c_1 y ^{-\beta} + o(y^{-\beta})$, the corresponding
endpoint behavior is
\[
  \widehat{q}(y(x)) = c_0 + c_1 \left(\frac{1-x}{2L}\right)^{\beta / \alpha} + o\bigl((1-x)^{\beta/\alpha}\bigr).
\]
If $\beta / \alpha$ is not an integer and $c_1 \neq 0$, this fractional
power prevents the mapped field from being differentiable at $x = 1$. The
same consideration applies to higher-order corrections: although smaller in
magnitude, they can still limit regularity through their higher-order
derivatives. Choosing $\alpha$ so that the correction exponents in the
reference variable $x$ are integers removes these fractional-power
obstructions and can substantially improve convergence. Section~\ref{sm-sec:exponents}
shows how $\alpha$ is chosen for each problem considered in this paper so that this condition holds.

The composite approximation~\eqref{eq:pou-blend} requires window functions
that are smooth, sum to one, and vanish outside their patches. We construct
them from a single cutoff profile. On the overlap $(y_a, y_b)$ of two
consecutive patches $K_{i-1}$ and $K_i$, we define
\begin{equation}\label{sm-eq:cutoff}
  \widetilde{\chi}(y;\, y_a, y_b) = \tfrac{1}{2}\,\operatorname{erfc}\!\bigl(6(2t-1)\bigr),
  \qquad t = \frac{y-y_a}{y_b-y_a},
\end{equation}
which is $C^\infty$ and decreases from one to zero, to within rounding,
across $[y_a,y_b]$~\cite{bruno2026two}. The two windows that share this
overlap are
\[
  \chi^{(i-1)}(y)
  = \widetilde{\chi}(y;\,y_a,y_b),
  \qquad
  \chi^{(i)}(y)
  = 1-\widetilde{\chi}(y;\,y_a,y_b), \qquad y \in [y_a, y_b].
\]
Elsewhere, each window equals one on the part of its patch that overlaps no
other patch and zero outside its patch. The overlap between the last local
patch and the asymptotic patch is treated in the same way. Since each window
is smooth, the composite approximation inherits the smoothness of the
patchwise interpolants.

\subsection{Two-dimensional representation}\label{sm-sec:2d}

Section~\ref{sec:methodology-general-2d} of the main text describes the
two-dimensional representation qualitatively: a chain of Cartesian patches
near the origin, a chain of polar patches in the far field, and a radial
partition of unity that joins the two. This section gives the corresponding
definitions. Near the origin, the local patches are the overlapping
rectangles
\begin{equation}\label{sm-eq:2d-local-patches}
  K^{(i)}_\mathrm{local}
  =
  \bigl(y_{1,\mathrm{start}}^{(i)},
        y_{1,\mathrm{end}}^{(i)}\bigr)
  \times (0,y_{2,\mathrm{end}}),
  \qquad i=1,\dots,m,
\end{equation}
which overlap only their immediate neighbors in $y_1$ and share a common
extent in $y_2$. On the closure of each patch, the restriction
$q^{(i)}=q|_{\overline{K^{(i)}_\mathrm{local}}}$ is approximated by the tensor-product
interpolant $q_{\mathbf n_i}^{(i)}$ on an $n_{i,1}\times n_{i,2}$
Chebyshev--Lobatto grid, where $\mathbf n_i=(n_{i,1},n_{i,2})$.

In the far field, mapping $y_1$ and $y_2$ separately to bounded intervals
can reduce the regularity of the normalized fields where both coordinates
tend to infinity. We therefore use polar coordinates $y_1=r\cos\theta$,
$y_2=r\sin\theta$ and cover the region $r>r_{\mathrm{start}}$ by a chain of
$m_{\mathrm{asym}}$ angular patches
\[
  K_{\mathrm{asym}}^{(j)}
  =
  \left\{
    (r\cos\theta,r\sin\theta):
    r>r_{\mathrm{start}},\;
    \theta\in
    \bigl(\theta_{\mathrm{start}}^{(j)},
          \theta_{\mathrm{end}}^{(j)}\bigr)
  \right\},
  \qquad j=1,\dots,m_{\mathrm{asym}},
\]
whose angular intervals cover $(0,\pi/2)$ and overlap only their immediate
neighbors. The radial coordinate of each patch is mapped to
$\rho\in[-1,1]$ by the two-dimensional analogue of~\eqref{sm-eq:map},
\begin{equation}\label{sm-eq:2d-map}
  \rho(r)
  =
  \frac{r^\alpha-r_{\mathrm{start}}^\alpha-L}
       {r^\alpha-r_{\mathrm{start}}^\alpha+L},
  \qquad
  L=r_{\mathrm{mid}}^\alpha-r_{\mathrm{start}}^\alpha,
\end{equation}
with $\alpha$ chosen as in Section~\ref{sm-sec:exponents}. With the
envelope $E_q(r)=r^\sigma$, the normalized field
\begin{equation}\label{sm-eq:2d-qhat}
  \widehat q^{(j)}(r,\theta)
  =
  \frac{q(r\cos\theta,r\sin\theta)}{E_q(r)},
  \qquad
  (r\cos\theta,r\sin\theta)
  \in\overline{K_{\mathrm{asym}}^{(j)}},
\end{equation}
is approximated on each patch by the tensor-product interpolant
$\widehat q_{\mathbf n_{\mathrm{asym}}^{(j)}}^{(j)}$ in $(\rho,\theta)$ on an
$n_\rho^{(j)}\times n_\theta^{(j)}$ Chebyshev--Lobatto grid on
$[-1,1]\times [\theta_{\mathrm{start}}^{(j)}, \theta_{\mathrm{end}}^{(j)}]$,
where $\mathbf n_{\mathrm{asym}}^{(j)}=(n_\rho^{(j)},n_\theta^{(j)})$.

As outlined in the main text, the two families of interpolants are combined
in two steps. Within each family, the interpolants are blended by
one-dimensional windows constructed from the cutoff~\eqref{sm-eq:cutoff}:
windows $\chi^{(i)}(y_1)$ in $y_1$ for the local patches and windows
$\chi_{\mathrm{polar}}^{(j)}(\theta)$ in $\theta$ for the asymptotic patches.
The two families are then joined by a radial transition. We choose a radius
$r_{\mathrm{blend}}$ satisfying
\[
  r_{\mathrm{start}}<r_{\mathrm{blend}}
  \leq
  \min\bigl(y_{1,\mathrm{end}}^{(m)},y_{2,\mathrm{end}}\bigr),
\]
so that every point of the quarter plane with
$r_{\mathrm{start}}<r<r_{\mathrm{blend}}$ lies in both the local and the
asymptotic region, and we construct from~\eqref{sm-eq:cutoff} a function
$\chi_r(r)$ that decreases from one to zero on
$[r_{\mathrm{start}},r_{\mathrm{blend}}]$. The two-dimensional windows are
\[
  \chi_{\mathrm{local}}^{(i)}(y_1,y_2)
  =
  \chi_r(r)\chi^{(i)}(y_1),
  \qquad
  \chi_{\mathrm{asym}}^{(j)}(y_1,y_2)
  =
  \bigl(1-\chi_r(r)\bigr)
  \chi_{\mathrm{polar}}^{(j)}(\theta),
\]
and the composite approximation is
\begin{equation}\label{sm-eq:2d-pou-blend}
  q_{\mathbf n}(y)
  =
  \sum_{i=1}^{m}
    \chi_{\mathrm{local}}^{(i)}(y)\,
    q_{\mathbf n_i}^{(i)}(y)
  +
  E_q(r)
  \sum_{j=1}^{m_{\mathrm{asym}}}
    \chi_{\mathrm{asym}}^{(j)}(y)\,
    \widehat q_{\mathbf n_{\mathrm{asym}}^{(j)}}^{(j)}
    (r,\theta),
\end{equation}
where $\mathbf n$ collects the resolutions of all patches.

\subsection{Loss and linear solvers}\label{sm-sec:loss}

The loss described in Section~\ref{sec:cheb-background} of the main text
combines collocation residuals of the governing equations, a penalty on
their spatial derivatives, and matching conditions on the patch overlaps.
We now specify each block of the residual vector for the one-dimensional
representation; the two-dimensional computation uses the same blocks on the
tensor-product grids, with the modifications given in
Section~\ref{sm-sec:boussinesq}. Let $f[\,\cdot\,](y)$ denote the residual of
a governing equation and $\kappa\geq1$ an integer refinement factor. The
equation blocks below are assembled for every governing equation, and the
matching block for every unknown field.
\begin{enumerate}
  \item \textbf{Local collocation residuals.}
  On each local patch $K_i$, the residual $f$ is evaluated at the $\kappa
  n_i$ Chebyshev--Lobatto nodes of that patch, and the block is scaled by
  $1/\sqrt{\kappa n_i}$.

  \item \textbf{Local gradient penalty.}
  On each local patch, the residual values of the preceding block are
  interpolated by a Chebyshev polynomial, whose derivative with respect to
  $y$ is evaluated at the same nodes. The block is scaled by
  $\sqrt{c_{\mathrm{fgrad}}}/\sqrt{\kappa n_i}$.

  \item \textbf{Asymptotic collocation residuals.}
  On the asymptotic patch, $f$ is evaluated at the nodes of the
  interpolation grid other than the node at $y=\infty$, and the block is
  scaled by $1/\sqrt{n_{\mathrm{asym}}}$. No gradient penalty is applied on
  this patch.

  \item \textbf{Patch-matching residuals.}
  For neighboring local patches $K_i$ and $K_{i+1}$, let
  $\{y_k^{\mathrm{ov}}\}_{k=1}^{n_{\mathrm{ov}}}$ denote the interpolation
  nodes of both patches that lie in $K_i\cap K_{i+1}$. The matching
  residuals are
  \[
    r_k^{\mathrm{match}}
    =
    \frac{
      q_{n_i}^{(i)}(y_k^{\mathrm{ov}})
      -
      q_{n_{i+1}}^{(i+1)}(y_k^{\mathrm{ov}})
    }{\sqrt{n_{\mathrm{ov}}}}, \qquad k=1, \dots, n_\mathrm{ov}.
  \]
  On the overlap of the last local patch with the asymptotic patch,
  $q_{n_m}^{(m)}$ is compared in the same way with the reconstructed field
  $E_q\,\widehat q_{n_{\mathrm{asym}}}$.
\end{enumerate}
The scaling of each block by the square root of its number of evaluation
points makes its contribution to the loss a mean square, so that patches
with different resolutions are weighted comparably.

Each Gauss--Newton iteration requires the solution of a damped linear
least-squares problem, with the Levenberg--Marquardt damping parameter
adapted automatically by \textsc{Optimistix}. For all but the largest
problems, we form the Jacobian of the residual vector by automatic
differentiation and solve each damped subproblem by dense QR
factorization~\cite{rader2023lineax}, which avoids forming the normal
equations. At the finest Boussinesq resolution and for the Euler
computation of Section~\ref{sm-sec:euler}, we instead avoid storing the
Jacobian and compute the step with LSMR~\cite{fong2011lsmr,paige1982lsqr},
which requires only products with the Jacobian and its transpose.
Unpreconditioned LSMR did not reach the required accuracy within a practical
number of iterations, and we therefore precondition it following
Blendenpik~\cite{avron2010blendenpik}. The preconditioner is the triangular
factor of the QR factorization of a row sketch $S = PFD$ of the Jacobian, in
which $D$ applies random sign flips, $F$ is a discrete cosine transform, and
$P$ samples $\lceil\gamma n\rceil$ rows uniformly at random, with $n$ the
number of unknowns. The sketch is assembled from transpose
Jacobian--vector products, so the Jacobian is never formed, and the
conditioning of the triangular factor is monitored by a $1$-norm
estimator~\cite{higham2000block}.

\subsection{Adaptive refinement}\label{sm-sec:adaptivity}

The adaptive procedure of Section~\ref{sec:general-adaptive} of the main
text uses the normalized tail maximum $e_k$ of the Chebyshev coefficients of
each local interpolant to decide whether a patch should be split, merged
with a neighbor, or reduced in degree. This section states the bound that
relates the coefficient tail to the interpolation error, and it lists the
settings of the procedure that are common to the Burgers and De Gregorio
computations.

For a field with an absolutely convergent Chebyshev
series~\eqref{sm-eq:cheb-coeffs}, the interpolant $q_n$ on $n$
Chebyshev--Lobatto nodes satisfies
\[
  q(y)-q_n(y)
  =
  \sum_{k=n}^{\infty}
    a_k\bigl(T_k(y)-T_{k^*}(y)\bigr),
\]
where $k^*\in\{0,\dots,n-1\}$ is the unique index with
$k\equiv\pm k^*\pmod{2(n-1)}$~\cite{trefethen2013approximation}. Since
$|T_k(y)|\leq1$ on $[-1,1]$, it follows that
\[
  \|q-q_n\|_\infty
  \leq 2\sum_{k=n}^{\infty}|a_k|,
\]
so the interpolation error is controlled by the tail of the coefficient
sequence. In practice, the criterion is applied to the coefficients of the
computed interpolant rather than to those of the unknown solution, and
their decay is therefore a practical diagnostic of resolution rather than a
rigorous error bound.

In the computations of Sections~\ref{sec:burgers} and~\ref{sec:degregorio},
the criterion is applied to the trailing coefficients: a patch with $n$
nodes is declared resolved when $e_{n-n_{\mathrm{tail}}}\le\varepsilon$, that
is, when its last $n_{\mathrm{tail}}$ normalized coefficients lie below
$\varepsilon$. The chop pass does not reduce a patch below $8$ nodes. Rather
than using a single tolerance, we run the adaptive cycle once for each
tolerance of a decreasing sequence. Throughout, the first local patch starts
at $y=0$, the extent $y^{(m)}_{\mathrm{end}}$ of the local region is held
fixed, and the asymptotic patch is not adapted. The tolerances and the
values of $N_{\max}$, $n_{\mathrm{tail}}$ and $\eta$ used for each problem
are listed in Tables~\ref{sm-tab:burgers-cfg} and~\ref{sm-tab:dg-cfg}.

\section{Parity extension}\label{sm-sec:parity-extension}

The fields are computed on the half-line $y\ge 0$ (or the quarter-plane
$y_1,y_2\ge 0$, or $r\ge 0$ for the Euler equations) and then extended
across the symmetry axes to obtain the full solution profiles for
validation and visualization. All residuals and errors reported in
Section~\ref{sec:results} of the main text, and the profiles shown in
Figures~\ref{fig:dg-profiles}, \ref{fig:bq-fields}
and~\ref{fig:bq-residuals}, are evaluated on these extended fields. This
section describes how the extensions are constructed to preserve smoothness
and numerical accuracy.

Since the computed fields satisfy the symmetry conditions only to within
the accuracy of the optimization, a direct reflection would in general
introduce a discontinuity in some derivative at the axis. We therefore blend
the computed field with its reflection over a narrow band around the axis.
Let $x$ denote the coordinate normal to the axis ($y$, $y_1$, or the
computational radius $r_c$ of Section~\ref{sm-sec:euler}), let $s=+1$ for an
even field and $s=-1$ for an odd one, and let $H$ denote the computed
composite approximation, evaluated on its innermost patch slightly beyond
the axis, $x\ge-\delta$. The extended field is
\begin{equation}\label{sm-eq:mirror}
  F(x) =
  \begin{cases}
    H(x), & x \ge \delta,\\
    w(x)\,H(x) + w(-x)\,s\,H(-x), & |x| < \delta,\\
    s\,H(-x), & x \le -\delta,
  \end{cases}
  \qquad
  w(x) = \frac{\widetilde\chi(-x;\,-\delta,\delta)}
              {\widetilde\chi(x;\,-\delta,\delta)+\widetilde\chi(-x;\,-\delta,\delta)},
\end{equation}
where $\widetilde\chi$ is the cutoff~\eqref{sm-eq:cutoff}. Since
$w(x)+w(-x)=1$, the extension satisfies $F(-x) = sF(x)$ exactly and
coincides with $H$ for $x\ge\delta$. Its first and second derivatives follow
from~\eqref{sm-eq:mirror} by the product rule, using the closed-form
derivatives of $w$, and satisfy the parity relations to within rounding
($\le 4.5\times10^{-16}$ in all cases). In two dimensions, $H$ and its
derivatives are evaluated at $(\pm x, \cdot)$, and the transverse
derivatives are blended in the same way. In practice, we take the blend half-width to be
the position of the first Chebyshev node off the axis on the innermost
patch, $\delta = y^{(1)}_{\mathrm{end}}\bigl(1-\cos(\pi/(n_1-1))\bigr)/2$.
Table~\ref{sm-tab:mirror} lists $\delta$ and the largest deviation of $F$
from the computed field on the half domain for every profile.

The $\lambda$-inference procedure of Remark~\ref{rem:burgers-lambda}
evaluates the residual near the origin and therefore uses a variant of this
construction that yields a single polynomial across the axis: the origin
patch is re-interpolated on
$[-y^{(1)}_{\mathrm{end}}, y^{(1)}_{\mathrm{end}}]$ from the
parity-reflected values of $H$, on a Chebyshev grid with twice as many
nodes, and the residual is evaluated on this interpolant.

\begin{table}[htb]
  \centering
  \caption{Blend half-width $\delta$ of~\eqref{sm-eq:mirror} and largest
    relative deviation $\max|F-H|/\max|H|$ of the extended field from the
    computed field on the half domain, over all fields.}
  \label{sm-tab:mirror}
  \begin{tabular}{lcc}
    \toprule
    Profile & $\delta$ & $\max|F-H|/\max|H|$ \\
    \midrule
    Burgers, $\lambda=1/2$, $1/4$, $1/6$
      & $2.6$, $1.4$, $1.1$ ($\times10^{-3}$) & $\le 6.9\times10^{-17}$ \\
    \texttt{clm} & $4.8\times10^{-3}$ & $5.7\times10^{-17}$ \\
    \texttt{ccf0} & $9.9\times10^{-4}$ & $5.8\times10^{-17}$ \\
    \texttt{ccf1} & $1.1\times10^{-2}$ & $7.0\times10^{-17}$ \\
    \texttt{ccf2} & $4.0\times10^{-3}$ & $2.9\times10^{-13}$ \\
    Boussinesq & $9.3\times10^{-3}$ & $1.3\times10^{-12}$ \\
    Euler & $1.9\times10^{-2}$ & $8.7\times10^{-8}$ \\
    \bottomrule
  \end{tabular}
\end{table}

\section{Far-field and mapping exponents}\label{sm-sec:exponents}

Before turning to the individual problems in
Sections~\ref{sm-sec:burgers}--\ref{sm-sec:boussinesq}, we derive the two
exponents that each of them must prescribe for every field, as listed at the
beginning of Section~\ref{sec:results} of the main text: the far-field
exponent $\sigma_q$, which defines the envelope $E_q = y^{\sigma_q}$
in~\eqref{eq:qhat}, and the mapping exponent $\alpha$, which is chosen to
keep the normalized fields smooth in the reference variable. Both exponents
follow from the asymptotic structure of the profile equations. We first
state the criteria that determine them and then apply these criteria to the
Burgers, De Gregorio and Boussinesq problems. The resulting values, which the
problem-specific sections below use without further comment, are collected
in Table~\ref{sm-tab:exponents}.

The exponent $\sigma_q$ is the leading power of $q$ as $y\to\infty$ and
follows from a dominant balance. In each profile equation considered here,
the terms that dominate at large $y$ form a linear combination
$c\,q + (1+\lambda)\,y\,\partial_y q$ acting on a single field, where the
constant $c$ depends on the equation. A power law $q \sim y^{\sigma}$
annihilates this combination if $\sigma = -c/(1+\lambda)$, and the exponents
of the remaining fields then follow from the relations that couple them to
the first. Only $\sigma_q$ is prescribed: the amplitude $\widehat q(\infty)$
and, in two dimensions, its angular dependence remain unknowns of the
optimization.

As shown in Section~\ref{sm-sec:grid}, a correction $y^{-\mu}$ in
$\widehat q$ appears in the reference variable as $(1-x)^{\gamma}$ with
$\gamma = \mu/\alpha$. When $\gamma$ is not an integer, this branch point at
$x=1$ limits the Chebyshev coefficients of $\widehat q$ to algebraic
decay~\cite{trefethen2013approximation},
\begin{equation}\label{sm-eq:branch-rate}
  |a_k| = O\bigl(k^{-(2\gamma+1)}\bigr).
\end{equation}
We therefore choose $\alpha$ so that as many correction exponents of
$\widehat q$ are integer multiples of $\alpha$ as possible, and we now identify these
exponents for each problem.

For the generalized De Gregorio and Boussinesq systems, factoring out the
envelopes leaves residuals whose terms differ from the leading balance by
integer powers of $y^{-1/(1+\lambda)}$, or of $r^{-1/(1+\lambda)}$ in two
dimensions. For example, substituting
$\Omega = y^{\sigma_\Omega}\widehat\Omega$ and
$U = y^{\sigma_U}\widehat U$ into the residual~\eqref{eq:dg-f1} gives
\begin{equation}\label{sm-eq:f1-normalized}
  y^{1/(1+\lambda)} f_1
  = (1+\lambda)\,y\,\partial_y\widehat\Omega
    + y^{-1/(1+\lambda)}\,N\bigl[\widehat\Omega, \widehat U\bigr],
\end{equation}
where $N$ is polynomial in $\widehat\Omega$, $\widehat U$ and their
scaled derivatives $y\,\partial_y \widehat\Omega$ and $y\,\partial_y \widehat U$. In the variable $s = y^{-1/(1+\lambda)}$,
equation~\eqref{sm-eq:f1-normalized} is an ordinary differential equation
whose right-hand side is analytic at $s=0$ and whose initial value
$\widehat\Omega(0)$ is free. Its solutions are therefore analytic in $s$,
and we take $\alpha = 1/(1+\lambda)$ for both problems.

For the Burgers equation, the normalized profile is known in closed form.
Substituting $U = y^{\sigma_U}\widehat U$ into~\eqref{eq:burgers-exact} with
$C=1$ gives $|\widehat U|^{1+1/\lambda} + y^{-1/(1+\lambda)}|\widehat U| = 1$,
so that $\widehat U$ is analytic in
$y^{-1/(1+\lambda)} = (y^{-\alpha})^{1/\lambda}$ with
$\alpha = \lambda/(1+\lambda)$. On the smooth
branches~\eqref{eq:burgers-smooth-lam}, $1/\lambda = 2i+2$ is an integer, so
this is an integer power of $y^{-\alpha}$ and the criterion is satisfied.

\begin{table}[htb]
  \centering
  \caption{Far-field exponents $\sigma_q$ and mapping exponents $\alpha$.
    For the generalized De Gregorio family the same $\alpha$ is used for
    both fields, and for the Boussinesq equations for all five fields.}
  \label{sm-tab:exponents}
  \begin{tabular}{lll}
    \toprule
    Problem & $\sigma_q$ & $\alpha$ \\
    \midrule
    Burgers & $\sigma_U = \lambda/(1+\lambda)$ & $\lambda/(1+\lambda)$ \\
    De Gregorio & $\sigma_\Omega = -1/(1+\lambda)$, \
      $\sigma_U = \lambda/(1+\lambda)$ & $1/(1+\lambda)$ \\
    Boussinesq & $\sigma_\Omega = -1/(1+\lambda)$, \
      $\sigma_\Phi = \sigma_\Psi = -2/(1+\lambda)$, \
      $\sigma_{U_1} = \sigma_{U_2} = \lambda/(1+\lambda)$ & $1/(1+\lambda)$ \\
    \bottomrule
  \end{tabular}
\end{table}

\section{The Burgers equation}\label{sm-sec:burgers}

Section~\ref{sec:burgers} of the main text reports the profiles for
$\lambda=1/2$, $1/4$ and $1/6$, together with their residuals and their
errors with respect to the exact solutions~\eqref{eq:burgers-exact}. This
section gives the discretizations on which these results were obtained. The
single field $U$ is represented on a chain of local patches covering
$[0, y^{(m)}_{\mathrm{end}}]$ and one asymptotic patch, with the far-field
and mapping exponents of Table~\ref{sm-tab:exponents}. The local patches are
produced by the adaptive procedure of Section~\ref{sec:general-adaptive},
starting from a single patch, with the settings described in
Section~\ref{sm-sec:adaptivity}. Table~\ref{sm-tab:burgers-cfg} lists the
configurations and the resulting discretizations; in each case the
procedure terminates with at most two local patches.

\begin{table}[htb]
  \centering
  \caption{Burgers configurations and final discretizations. Node counts
    are numbers of Chebyshev--Lobatto points. The lower block is common to
    all three profiles.}
  \label{sm-tab:burgers-cfg}
  \begin{tabular}{lccc}
    \toprule
    $\lambda$ & $1/2$ & $1/4$ & $1/6$ \\
    \midrule
    $\alpha = \sigma_U$ & $1/3$ & $1/5$ & $1/7$ \\
    local patches $[y_{\mathrm{start}}, y_{\mathrm{end}}]$
      & $[0,3]$
      & $[0,1.667]$, $[1.333,3]$
      & $[0,0.926]$, $[0.741,3]$ \\
    local node counts & $54$ & $56$, $23$ & $47$, $42$ \\
    $n_{\mathrm{asym}}$ & $32$ & $32$ & $40$ \\
    \midrule
    \multicolumn{4}{l}{$y^{(m)}_{\mathrm{end}} = 3$,
      $y^{\mathrm{asym}}_{\mathrm{start}} = 2$,
      $y^{\mathrm{asym}}_{\mathrm{mid}} = 10$; $\kappa = 3$,
      $c_{\mathrm{fgrad}} = 1$} \\
    \multicolumn{4}{l}{$\varepsilon \in \{10^{-14},\,10^{-15},\,6\times10^{-16}\}$;
      $N_{\max} = 64$, $n_{\mathrm{tail}} = 5$, $\eta = 0.2$;
      at most $40$ adaptive rounds} \\
    \multicolumn{4}{l}{initial grid: one patch $[0,3]$ with $64$ nodes;
      $100$ \ac{LM} iterations per round} \\
    \multicolumn{4}{l}{final solve: at most $2000$ \ac{LM} iterations,
      absolute tolerance $10^{-16}$} \\
    \bottomrule
  \end{tabular}
\end{table}

\subsection{Comparison with PINNs}\label{sm-sec:pinns}
Our \ac{PINN} implementation follows the methodology of~\cite{wang2025discovery}, including a solution ansatz that incorporates symmetry and far-field decay through mapped input coordinates $(q, \beta)$, fully connected neural networks with $\tanh$ activations, residual and residual-derivative losses evaluated at both prescribed and adaptively sampled collocation points, and a two-stage training procedure with Fourier features in the second stage. For optimization, we use the full-matrix stochastic Gauss--Newton optimizer implemented in \texttt{kfac-jax}~\cite{kfac-jax2022github}, with automatic adaptation of the learning rate, momentum, and damping parameters. Since not all implementation details are specified in~\cite{wang2025discovery}, we describe our choices and deviations below.

Both training stages use fully connected networks with $\tanh$ hidden activations and linear output layers, rather than exponential-type output activations. In the first stage, the mapped coordinate $q$ is passed through a multilayer perceptron with hidden width $w$ to obtain $N_1(q)$. In the second stage, $q$ is mapped to $m=w$ fixed random Fourier frequencies through the features $[\cos(qB),\sin(qB)]$, where $B\in\mathbb{R}^{1\times m}$ has independent entries sampled from $\mathcal{N}(0,\sigma^2)$. These features are passed through two hidden layers of width $w$ and a linear output layer to obtain $N_2(q)$. The Fourier matrix $B$ is held fixed during training. For $\lambda=1/2$, we consider network widths $w=8,16,32$, while for $\lambda=1/4$ and $1/6$, we use $w=16,24,32$. The first-stage networks have two hidden layers for widths $8$, $16$, and $24$, and three hidden layers for width $32$. Across these architectures, the first-stage networks contain $97$, $321$, $673$, or $2{,}209$ trainable parameters, respectively, while the corresponding second-stage networks contain $217$, $817$, $1{,}801$, or $3{,}169$ trainable parameters.

The loss includes the equation residual weighted by $(1+y^2)^{-\alpha/2}$, corresponding to its normalized far-field form. We also penalize the derivative of the residual with respect to $z=\operatorname{arcsinh}(y)$, evaluated at the same collocation points as the residual itself. No second-derivative loss is included. Each training batch contains $3{,}800$ collocation points on $y\geq0$, comprising $1{,}500$ points sampled uniformly in $z\in[0,40]$, $1{,}500$ sampled uniformly in $y\in[0,10]$, and $300$ sampled uniformly in $y\in[0,0.1]$. The remaining $500$ points are selected adaptively from a pool of $8{,}000$ points sampled uniformly in $z\in[0,40]$, with selection probabilities proportional to the residual magnitude. These selected points are then perturbed in $z$ by independent Gaussian noise with standard deviation $0.01$. The entire batch is resampled every $2{,}000$ training iterations.

The second training stage introduces a correction of the form $\varepsilon N_2$, where $\varepsilon$ is determined by the maximum residual from the first stage on a fixed probe grid. The output layer of $N_2$ is initialized to zero, and the second-stage loss is scaled by $1/\varepsilon$. We use a fixed Fourier-feature scale of $\sigma=2\pi$, rather than selecting it from the dominant frequency of the first-stage residual. The latter approach produced a scale of approximately $\sigma=150$ for $\lambda=1/4$, for which the second-stage network did not train successfully.

For both stages, all trainable parameters are treated as a single full curvature block in the stochastic Gauss--Newton optimizer. The curvature estimate uses an exponential-moving-average decay of $0.999$ and is inverted at every iteration. The damping parameter is initialized at $10^{-3}$. If the residual on the probe grid diverges, the training procedure restores the best previously obtained parameters and restarts the optimizer state, although this safeguard was not triggered in the reported experiments.

The first stage is trained for $30{,}000$ iterations. The second stage terminates when the maximum residual on the probe grid improves by less than $5\%$ over $20{,}000$ iterations, subject to a maximum of $180{,}000$ iterations. The parameters corresponding to the smallest observed probe residual are retained. The probe grid contains points distributed across both the mapped far-field region and the local region, allowing the stopping criterion to account for errors at different spatial scales.

\begin{table}[t]
\centering
\caption{Comparison of \ac{ODIL} and \acp{PINN} for the Burgers
self-similar profiles across different neural network sizes.
Maximum residuals and relative maximum errors are evaluated
over the ``infinite domain'' validation region.}
\label{tab:burgers-comparison-full}

\small
\renewcommand{\arraystretch}{1.3}
\setlength{\tabcolsep}{6pt}

\begin{tabular}{@{}c l r r r r@{}}
\toprule
$\lambda$
& Method
& Parameters (Stage 1/2)
& Time (s)
& Max. $f$
& Max. $U - U_\mathrm{ex}$ \\
\midrule

$1/2$
& ODIL
& 84/--
& 0.405
& $5.64\times10^{-14}$
& $1.23\times10^{-15}$ \\

& PINN
& 97/217
& 295
& $7.14\times10^{-14}$
& $4.25\times10^{-15}$ \\

& PINN
& 321/817
& 159
& $9.22\times10^{-15}$
& $4.52\times10^{-15}$ \\

& PINN
& 2,209/3,169
& 529
& $9.16\times10^{-15}$
& $3.61\times10^{-15}$ \\

\midrule

$1/4$
& ODIL
& 109/--
& 0.620
& $1.60\times10^{-14}$
& $8.96\times10^{-16}$ \\

& PINN
& 321/817
& 572
& $1.02\times10^{-13}$
& $6.03\times10^{-15}$ \\

& PINN
& 673/1,801
& 504
& $5.46\times10^{-14}$
& $2.81\times10^{-15}$ \\

& PINN
& 2,209/3,169
& 1,430
& $3.04\times10^{-14}$
& $2.01\times10^{-15}$ \\

\midrule

$1/6$
& ODIL
& 127/--
& 0.754
& $1.62\times10^{-14}$
& $7.59\times10^{-16}$ \\

& PINN
& 321/817
& 486
& $1.83\times10^{-13}$
& $1.05\times10^{-14}$ \\

& PINN
& 673/1,801
& 872
& $8.93\times10^{-14}$
& $5.40\times10^{-15}$ \\

& PINN
& 2,209/3,169
& 1,360
& $4.53\times10^{-14}$
& $2.84\times10^{-15}$ \\

\bottomrule
\end{tabular}%

\end{table}

\section{The generalized De Gregorio family}\label{sm-sec:degregorio}

Section~\ref{sec:degregorio} of the main text computes the \ac{CLM} profile
and three branches of the \ac{CCF} model by a two-stage optimization, with
the Hilbert transform in the residual~\eqref{eq:dg-f2} evaluated by
quadrature. This section describes the normalization imposed in each stage,
lists the configurations and final discretizations, and details the
quadrature.

\subsection{Normalization and two-stage solve}\label{sm-sec:dg-stages}

The system~\eqref{eq:dg-f1}--\eqref{eq:dg-f2} is invariant under the
rescaling $\Omega(y)\mapsto\Omega(\mu y)$, $U(y)\mapsto\mu^{-1}U(\mu y)$,
$\mu>0$, under which the far-field amplitude transforms as
$\widehat\Omega(\infty)\mapsto\mu^{-1/(1+\lambda)}\widehat\Omega(\infty)$.
A single scalar condition is therefore required to select one member of
each family of solutions. In stage~1, this condition is imposed as a
penalty weighted by $\sqrt{c_g}$: we prescribe $\partial_y\Omega(0) = -2$
for \texttt{clm} and \texttt{ccf0}, consistent with the exact \ac{CLM}
profile~\eqref{eq:clm-exact}, and $\Omega(0.5) = -0.05$ for \texttt{ccf1}.
Stage~2 replaces the penalty by a constraint on the far field:
$\widehat\Omega(\infty)$ is fixed at its stage-1 value and removed from the
optimization variables, the penalty is dropped, and $c_{\mathrm{fgrad}}$ is
set to zero. The \texttt{ccf2} solve differs in that it is initialized from
a separately computed profile and uses no stage-1 penalty, with
$c_{\mathrm{fgrad}} = 1$; its scaling is fixed only in stage~2. In all
cases, the local patches are produced by the adaptive procedure with the
settings of Section~\ref{sm-sec:adaptivity}, with $\Omega$ and $U$ sharing
their patch boundaries. Table~\ref{sm-tab:dg-cfg} lists the settings, and
Table~\ref{sm-tab:dg-patches} the resulting local patches.

\begin{table}[htb]
  \centering
  \small
  \setlength{\tabcolsep}{4pt}
  \caption{Configurations for the generalized De Gregorio family. The
    asymptotic grids of $\Omega$ and $U$ share $y^{\mathrm{asym}}_{\mathrm{start}}$,
    $y^{\mathrm{asym}}_{\mathrm{mid}}$ and $n_{\mathrm{asym}}$. All cases use
    $\kappa = 3$ and at most $10$ merge passes.}
  \label{sm-tab:dg-cfg}
  \begin{tabular}{lcccc}
    \toprule
    & \texttt{clm} & \texttt{ccf0} & \texttt{ccf1} & \texttt{ccf2} \\
    \midrule
    $a$ & $0$ & $-1$ & $-1$ & $-1$ \\
    $\lambda$ & $0$ & $1.180777662899$ & $0.6057337012032$ & $0.471324863862$ \\
    $\alpha$ & $1$ & $0.4586$ & $0.6228$ & $0.6797$ \\
    $y^{(m)}_{\mathrm{end}}$ & $6$ & $6$ & $11$ & $45$ \\
    $y^{\mathrm{asym}}_{\mathrm{start}}$, $y^{\mathrm{asym}}_{\mathrm{mid}}$
      & $5.5$, $50$ & $5.5$, $50$ & $10$, $20$ & $40$, $45$ \\
    $n_{\mathrm{asym}}$ & $40$ & $26$ & $80$ & $108$ \\
    \midrule
    stage-1 condition & $\partial_y\Omega(0)=-2$ & $\partial_y\Omega(0)=-2$
      & $\Omega(0.5)=-0.05$ & none \\
    stage-1 $c_{\mathrm{fgrad}} = c_g$ & $100$ & $100$ & $100$ & $1$ \\
    stage-1 / 2 \ac{LM} iterations & $50$ / $30$ & $50$ / $30$
      & $500$ / $50$ & $2000$ / $30$ \\
    \midrule
    $\log_{10}\varepsilon$ sequence & $-11,-13,-14$
      & $-11,-13,-14,-15$
      & $-11,-12,\dots,-15$ & $-8,-10,-11,-13$ \\
    $N_{\max}$, $n_{\mathrm{tail}}$, $\eta$ & $32$, $5$, $0.2$
      & $32$, $10$, $0.2$ & $64$, $10$, $0.2$ & $64$, $10$, $0.1$ \\
    rounds; \ac{LM} iterations per round & $40$; $100$ & $40$; $100$
      & $40$; $100$ & $12$; $300$ \\
    \bottomrule
  \end{tabular}
\end{table}

\begin{table}[htb]
  \centering
  \small
  \caption{Final local patches $[y_{\mathrm{start}}, y_{\mathrm{end}}]$ and
    node counts $n_\Omega$, $n_U$ for the generalized De Gregorio family.}
  \label{sm-tab:dg-patches}
  \begin{tabular}{lccc@{\qquad}lccc}
    \toprule
    & $[y_{\mathrm{start}}, y_{\mathrm{end}}]$ & $n_\Omega$ & $n_U$
    & & $[y_{\mathrm{start}}, y_{\mathrm{end}}]$ & $n_\Omega$ & $n_U$ \\
    \midrule
    \texttt{clm}  & $[0, 1.029]$      & 24 & 22 & \texttt{ccf2} & $[0, 3.168]$     & 45 & 30 \\
                  & $[0.823, 1.852]$  & 19 & 18 & & $[3.051, 4.220]$  & 54 & 43 \\
                  & $[1.481, 3.992]$  & 23 & 21 & & $[4.180, 4.255]$  & 19 & 15 \\
                  & $[3.490, 6]$      & 18 & 15 & & $[4.247, 4.322]$  & 42 & 32 \\
    \cmidrule(r){1-4}
    \texttt{ccf0} & $[0, 0.176]$      & 22 & 15 & & $[4.308, 4.450]$  & 47 & 39 \\
                  & $[0.141, 0.572]$  & 23 & 21 & & $[4.443, 5.407]$  & 37 & 33 \\
                  & $[0.457, 1.232]$  & 23 & 21 & & $[5.311, 11.758]$ & 54 & 32 \\
                  & $[1.077, 1.852]$  & 16 & 16 & & $[11.191, 16.864]$ & 44 & 38 \\
                  & $[1.481, 2.876]$  & 20 & 17 & & $[15.786, 26.566]$ & 38 & 34 \\
                  & $[2.597, 6]$      & 23 & 20 & & $[24.518, 35.298]$ & 39 & 14 \\
    \cmidrule(r){1-4}
    \texttt{ccf1} & $[0, 1.886]$      & 22 & 19 & & $[34.220, 39.894]$ & 50 & 39 \\
                  & $[1.509, 3.764]$  & 51 & 43 & & $[39.326, 45]$    & 64 & 64 \\
                  & $[3.554, 4.602]$  & 33 & 29 & & & & \\
                  & $[4.225, 11]$     & 54 & 44 & & & & \\
    \bottomrule
  \end{tabular}
\end{table}

\subsection{The Hilbert transform}\label{sm-sec:dg-ht}

The residual~\eqref{eq:dg-f2} requires the Hilbert transform of $\Omega$ at
every collocation and validation point. Because the transform is nonlocal
and $\Omega$ decays only algebraically, its evaluation requires a quadrature
that resolves both the near field and the tail and that remains accurate
near the singularity of the kernel. We first rewrite the transform in a
form with a removable singularity. For odd $\Omega$, folding
the integral~\eqref{eq:hilbert-def} about the origin and subtracting
$2y_0\Omega(y_0)/(y_0^2-y^2)$, whose principal value vanishes, gives
\begin{equation}\label{sm-eq:ht-kernel}
  \Hop[\Omega](y_0)
  = \frac{1}{\pi}\int_0^\infty K(y_0, y)\,dy,
  \qquad
  K(y_0,y) = \frac{2\bigl(y\,\Omega(y) - y_0\,\Omega(y_0)\bigr)}{y_0^{2}-y^{2}},
\end{equation}
whose kernel has the finite diagonal limit
\begin{equation}\label{sm-eq:ht-diagonal}
  K(y_0,y_0) = -\frac{\Omega(y_0) + y_0\,\partial_y\Omega(y_0)}{y_0}.
\end{equation}
At the origin, $\Hop[\Omega](0) = -(2/\pi)\int_0^\infty \Omega(y)/y\,dy$ is
evaluated by a separate operator, and the self-term of a quadrature node at
$y = 0$ uses the limit $K(0,0) = -2\,\partial_y\Omega(0)$. The diagonal
limit is substituted on a branch whose inactive numerator and denominator
are replaced by dummy values, so that automatic differentiation never
encounters $0/0$.

The integral~\eqref{sm-eq:ht-kernel} is evaluated by Clenshaw--Curtis
quadrature~\cite{trefethen2008gauss} on a grid that is independent of the
solution grid. It consists of contiguous, nonoverlapping local segments
covering $[0, y^{\mathrm{ht}}_{\mathrm{start}}]$ and one segment on
$[y^{\mathrm{ht}}_{\mathrm{start}},\infty]$, mapped by~\eqref{sm-eq:map} with
the same $\alpha$ as the asymptotic patch. The field $\Omega$ is sampled at
the quadrature nodes through its composite approximation, and its
derivative on each segment is obtained in that segment's Chebyshev basis.
Since the accuracy of the quadrature depends on the position of the target
$y_0$ relative to the segments, the targets are grouped by magnitude into
panels, each with its own quadrature grid placed so that $y_0$ lies well
inside the resolved region (Table~\ref{sm-tab:dg-ht}). When a target lies
near but not on a node, the difference $y - y_0$ suffers from cancellation;
to avoid this, the local contribution is computed at each segment's own
nodes and barycentrically interpolated to the targets.

On the mapped segment, the slowly decaying part of the kernel is treated
analytically. The term $2y_0\Omega(y_0)/y^2$ is subtracted from $K$ and
integrated in closed form, contributing
$2y_0\Omega(y_0)/y^{\mathrm{ht}}_{\mathrm{start}}$, and the remaining
integrand decays like $y^{-1-\alpha}$. Its value at the endpoint $x=1$
($y=\infty$) is $-\widehat\Omega(\infty)/(L\alpha)$, independent of $y_0$,
and $+\widehat\Omega(\infty)/(2L\alpha)$ for the operator at the origin,
with $L$ as in~\eqref{sm-eq:map}. Both values are evaluated with the current
$\widehat\Omega(\infty)$, so that the contribution of the tail is
differentiated with respect to the far-field amplitude. This construction
requires $\alpha\le1$, which holds for all cases in
Table~\ref{sm-tab:dg-cfg}.

We assess the accuracy of the quadrature by rebuilding every panel with all
node counts doubled and comparing the two evaluations of $\Hop[\Omega]$ at
the collocation targets. The residuals in
Table~\ref{tab:clm-ccf-results} of the main text are computed with the
validation panels of Table~\ref{sm-tab:dg-ht}, excluding points close to
quadrature nodes as described in Section~\ref{sec:degregorio}.

\begin{table}[htb]
  \centering
  \small
  \setlength{\tabcolsep}{4pt}
  \caption{Hilbert-transform quadrature panels: target range of $|y_0|$,
    local segments with node counts, and the mapped segment
    $(y^{\mathrm{ht}}_{\mathrm{start}}, y^{\mathrm{ht}}_{\mathrm{mid}})$.
    Local segments have $256$ nodes and mapped segments $512$, except for
    \texttt{ccf2}, whose local segments are
    $[0,4]$, $[4,4.18]$, $[4.18,4.45]$, $[4.45,4.8]$, $[4.8,y^{\mathrm{ht}}_{\mathrm{start}}]$
    with $512$, $64$, $512$, $96$, $224$ nodes and whose mapped segments
    have $768$ nodes. The validation panels are used for the residuals
    reported in the main text.}
  \label{sm-tab:dg-ht}
  \begin{tabular}{llll}
    \toprule
    Case & Stage & $|y_0|$ & local / mapped segment \\
    \midrule
    \texttt{clm}, \texttt{ccf0}
      & 1 & all & $[0,10]$ / $(10, 10^3)$ \\
      & 2 & $<10^5$; $\ge10^5$ & $[0,10]$ / $(10,10^3)$; \ $[0,10^4]$ / $(10^4,10^6)$ \\
    \texttt{ccf1}
      & 1 & all & $[0,10]$ / $(10, 10^3)$ \\
      & 2 & $<10^2$; $[10^2,10^5)$; $\ge10^5$
        & $[0,1]$ / $(1,10)$; \ $[0,10]$ / $(10,10^3)$; \ $[0,10^2]$ / $(10^2,10^4)$ \\
    \texttt{clm}, \texttt{ccf0}, \texttt{ccf1}
      & val. & $<10^4$; $[10^4,10^6)$ & $[0,10]$ / $(10,10^3)$; \ $[0,10^2]$ / $(10^2,10^4)$ \\
      & & $[10^6,10^{10})$; $\ge10^{10}$ & $[0,10^4]$ / $(10^4,10^6)$; \ $[0,10^9]$ / $(10^9,10^{11})$ \\
    \texttt{ccf2}
      & 1, 2 & $<10^5$; $\ge10^5$ & $y^{\mathrm{ht}}_{\mathrm{start}}=40$, $(40,10^3)$; \ $100$, $(100,10^3)$ \\
      & val. & $<10^4$; $[10^4,10^6)$ & $(40,10^3)$; \ $(100,10^5)$ \\
      & & $[10^6,10^{10})$; $\ge10^{10}$ & $(10^4,10^7)$; \ $(10^9,10^{12})$ \\
    \bottomrule
  \end{tabular}
\end{table}

\section{The two-dimensional Boussinesq equations}\label{sm-sec:boussinesq}

Section~\ref{sec:boussinesq} of the main text solves the Boussinesq profile
equations in vorticity form for five fields, using one local and one
asymptotic patch per field and no adaptive refinement. This section
discusses the formulation, specifies the discretization, and describes the
continuation from a coarse to a fine resolution by which the reported
profile was obtained.

\subsection{Formulation}\label{sm-sec:bq-formulation}

In the formulation~\eqref{eq:bq-residuals}, the five
fields~\eqref{eq:bq-unknowns} are independent unknowns. The temperature
$\Theta$ itself is not represented; only its gradient $(\Phi,\Psi)$ enters
the equations, and the residual $f_6$ constrains this pair to be a gradient
field. Likewise, the velocity components are primitive unknowns, so that the
formulation requires no stream function, Poisson solve or Biot--Savart
evaluation: the divergence-free condition $f_4$ and the definition of the
vorticity $f_5$ are imposed as residuals of the same least-squares problem
as the transport equations. Posing the problem on the quarter
plane~\eqref{eq:quarter-plane} places the origin on the boundary of the
computational domain, where the Chebyshev nodes cluster along both the
symmetry axis and the wall. On the asymptotic patch, the Neumann condition
for $\Psi$ and $U_2$ is applied to the normalized field $\widehat q$; since
$\partial_{y_1}r^{\sigma}$ vanishes on $y_1=0$, this is equivalent to the
condition on $q$.

\subsection{Discretization}\label{sm-sec:bq-disc}

Each field is represented on one Cartesian local patch $[0,15]^2$ and one
polar asymptotic patch covering $\theta\in[0,\pi/2]$, where $\theta=0$
corresponds to the wall and $\theta=\pi/2$ to the symmetry axis. The
asymptotic patch starts at $r_{\mathrm{start}}=10$ and uses the radial
map~\eqref{sm-eq:2d-map} with $r_{\mathrm{mid}}=14$ and the exponents of
Table~\ref{sm-tab:exponents}, and the two patches are blended radially on
$10\le r\le15$. The nodes at $\rho=1$ are unknowns of the optimization and
represent the angular profile $\widehat q(\infty,\theta)$ of each field at
infinity. The residuals are collocated on a grid refined by $\kappa=2$ in
every direction on both patches and shared by all six residuals; the nodes
at $\rho=1$ are excluded from this grid and from the patch-matching term.
Since the scaling exponent enters the transport
field~\eqref{eq:transport-field}, the far-field exponents and the mapping
exponent, any change of $\lambda$ requires the grids to be rebuilt.

\subsection{Optimization}\label{sm-sec:bq-optimize}

The profile reported in Section~\ref{sec:boussinesq} is computed in two
steps, on the coarse and fine discretizations of
Table~\ref{sm-tab:bq-cfg}. The coarse problem is solved from the zero
initial guess with dense QR factorization. Its solution is then transferred
to the fine grids by evaluating the composite approximation at the new
nodes and dividing by the envelopes. At the nodes at $r=\infty$, where the
ratio of envelopes is indeterminate in floating point, the old normalized
asymptotic interpolant is instead evaluated directly in $(\rho,\theta)$.
The transferred solution serves as the initial guess for the fine problem,
which is solved with the sketch-preconditioned LSMR solver of
Section~\ref{sm-sec:loss}.

\begin{table}[htb]
  \centering
  \caption{Boussinesq discretizations. Node counts are
    $n_1\times n_2$ on the local patch and $n_\rho\times n_\theta$ on the
    asymptotic patch. Both use $\lambda = 1.9205599746927$, the local
    patch $[0,15]^2$, $r_{\mathrm{start}} = 10$, $r_{\mathrm{mid}} = 14$,
    $\kappa = 2$, $c_{\mathrm{fgrad}} = c_g = c_{\mathrm{bc}} = 1$, and at
    most $500$ \ac{LM} iterations with absolute tolerance $10^{-14}$. The
    results of the main text use the fine discretization.}
  \label{sm-tab:bq-cfg}
  \begin{tabular}{lcc}
    \toprule
    & coarse & fine \\
    \midrule
    degrees of freedom & $8\,304$ & $36\,889$ \\
    local $n_1\times n_2$: $\Omega$, $\Psi$ & $40\times24$ & $64\times35$ \\
    \qquad $\Phi$ & $40\times24$ & $64\times40$ \\
    \qquad $U_1$, $U_2$ & $32\times32$ & $52\times40$ \\
    asymptotic $n_\rho\times n_\theta$: $\Omega$, $\Phi$ & $32\times40$ & $64\times88$ \\
    \qquad $\Psi$ & $16\times24$ & $64\times88$ \\
    \qquad $U_1$ & $16\times16$ & $64\times72$ \\
    \qquad $U_2$ & $16\times24$ & $64\times72$ \\
    \midrule
    initial guess & $q \equiv 0$ & reprojected coarse solution \\
    linear solver & dense QR & sketched LSMR, $\gamma = 1.2$ \\
    \qquad tolerances, max.\ iterations & --- & $10^{-13}$ rel., $10^{-15}$ abs., $3000$ \\
    \bottomrule
  \end{tabular}
\end{table}

\section{The three-dimensional axisymmetric Euler equations}\label{sm-sec:euler}

Sections~\ref{sec:boussinesq} and~\ref{sec:conclusions} of the main text
identify singularity formation in the unforced three-dimensional Euler
equations~\cite{ganeshram2026stable} as a target for extensions of the
present method. As a preliminary step in this direction, we compute an
axisymmetric self-similar profile with the Chebyshev \ac{ODIL} framework.
Unlike the computations of the main text, this computation uses no
asymptotic patch: the domain is instead mapped and truncated, and the
far-field behavior is imposed as a penalty.

\paragraph{Equations.} We use the variables $U = u^\theta/r$,
$\Omega = \omega^\theta/r$ and
$\Psi = \psi^\theta/r$~\cite{luo2014potentially}. For a prescribed scaling
exponent $\lambda$ and an unknown axial translation constant $C$, the
self-similar profile equations read
\begin{equation}\label{sm-eq:euler-residuals}
  \begin{aligned}
    f_1 &= U + a_r\,\partial_r U + a_z\,\partial_z U - 2U\,\partial_z\Psi,\\
    f_2 &= (1+\lambda)\,\Omega + a_r\,\partial_r\Omega + a_z\,\partial_z\Omega - 2U\,\partial_z U,\\
    f_3 &= \Omega + \bigl(\partial_{rr} + \tfrac{3}{r}\,\partial_r + \partial_{zz}\bigr)\Psi,
  \end{aligned}
  \qquad
  \begin{aligned}
    a_r &= \lambda r - r\,\partial_z\Psi,\\
    a_z &= C + \lambda z + 2\Psi + r\,\partial_r\Psi,
  \end{aligned}
\end{equation}
on $r\ge0$, $z\in\mathbb{R}$, with all fields even in $r$. We take
$\lambda = 1/2$. Linearizing $f_1$ and $f_2$ about zero and inverting the
operator in $f_3$ gives the far-field decay $U\sim R^{-1/\lambda}$,
$\Omega\sim R^{-1-1/\lambda}$ and $\Psi\sim R^{1-1/\lambda}$, where
$R=(r^2+z^2)^{1/2}$; at $\lambda=1/2$ these are $R^{-2}$, $R^{-3}$ and
$R^{-1}$, respectively.

\paragraph{Discretization.} The physical coordinates are mapped by
$r = \sinh(r_c)$ and $z = \sinh(z_c)$ to the computational box
$(r_c,z_c)\in[0,30]\times[-30,30]$, which truncates the domain at
$\sinh(30)\approx5.3\times10^{12}$. Each field is represented on two
overlapping tensor-product Chebyshev patches, which cover
$z \in [-\sinh 30, 2]$ and $[-2, \sinh 30]$ and the full radial range, each
with $64\times64$ nodes, and which are blended in $z_c$ by the partition of
unity of Section~\ref{sm-sec:grid}. All nodal values, including those on
the outer edges, are unknowns. Derivatives are computed in $(r_c,z_c)$ and
converted to physical coordinates through the map; on the axis, the term
$(3/r)\partial_r\Psi$ is replaced by its regular limit. The residuals are
collocated on a grid refined by $\kappa=2$.

\paragraph{Loss.} The residual vector consists of the following terms:
\begin{itemize}
  \item the residuals $f_1$, $f_2$ and $f_3$ on the residual grid and the
        gradients of their interpolants, weighted by
        $\sqrt{c_{\mathrm{fgrad}}}$, as in Section~\ref{sm-sec:loss};
  \item the patch-matching term on the overlap;
  \item the even-parity condition $\partial_{r_c}q = 0$ on $r_c=0$ for
        each field, weighted by $\sqrt{c_{\mathrm{axis}}}$;
  \item the far-field condition
        $r\,\partial_r q + z\,\partial_z q + s_q\,q = 0$ with
        $(s_U, s_\Omega, s_\Psi) = (2,3,1)$, consistent with the decay rates
        above, imposed at the residual nodes with $R > R_{\mathrm{ff}}$;
  \item the normalization $\int U\,dr_c\,dz_c = u_{\mathrm{int}}$,
        evaluated by Clenshaw--Curtis quadrature weighted by the partition
        of unity and weighted by $\sqrt{c_{\mathrm{norm}}}$.
\end{itemize}
The unknowns are the nodal values of $U$, $\Omega$ and $\Psi$ together with
the constant $C$. The normalization fixes the scaling freedom of the
profile equations and makes the zero field a nonstationary point of the
loss, so that the solve can be started from $q\equiv0$ and $C=0$. We use
$c_{\mathrm{fgrad}} = c_{\mathrm{axis}} = c_{\mathrm{norm}} = 1$, unit weight
for the far-field term, $R_{\mathrm{ff}} = 3.6\times10^{10}$ and
$u_{\mathrm{int}} = 1000$. The problem is solved by at most $1000$ \ac{LM}
iterations with the sketched LSMR solver of Section~\ref{sm-sec:loss}
($\gamma = 1.2$, tolerances $10^{-13}$ relative and $10^{-15}$ absolute, at
most $3000$ iterations), and the converged translation constant is
$C = -6.609$.

\begin{figure}[htb]
  \centering
  \includegraphics[width=0.85\textwidth]{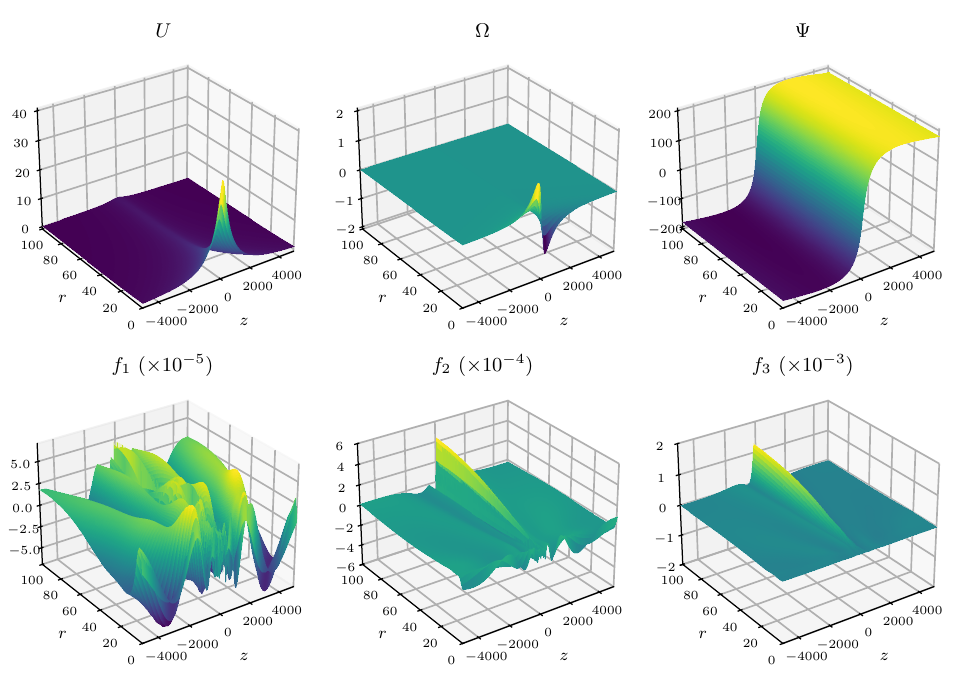}
  \caption{Fields (top) and residuals (bottom) of the axisymmetric Euler
    profile on $0\le r\le100$, $|z|\le5000$. Each residual is scaled by
    the power of ten given in its panel title.}
  \label{sm-fig:euler}
\end{figure}

\paragraph{Results.} The residuals are evaluated on the even extension of
Section~\ref{sm-sec:parity-extension} at $10^6$ points in each of four
regions, sampled uniformly in $(r_c,z_c)$: the full box $[-30,30]^2$, the
disk $r_c^2+z_c^2<9$, the axis $r_c=0$, and the far field
$R > R_{\mathrm{ff}}$. Table~\ref{sm-tab:euler-results} reports the results,
and Figure~\ref{sm-fig:euler} shows the fields and residuals near the core.
The relative tails of its Chebyshev coefficients reach
$1.6\times10^{-4}$ for $\Omega$ and $3.8\times10^{-5}$ for $\Psi$. The
discretization, whose patches are split at $z=\pm2$, is therefore not
matched to the profile, and the solution is not spectrally converged.
Beyond $|z|\approx10^5$, the computed fields are dominated by noise of order
$3\times10^{-6}$.
These limitations point to the extensions discussed in
Section~\ref{sec:conclusions}: adaptive refinement in two dimensions, which
would align the patches with the elongated core, and an asymptotic patch in
place of the truncated domain. Despite these limitations, the computed profiles already attain residuals comparable to or better than those reported in~\cite{ganeshram2026stable}.
\begin{table}[htb]
  \centering
  \small
  \setlength{\tabcolsep}{4pt}
  \caption{Residuals~\eqref{sm-eq:euler-residuals} of the axisymmetric Euler
    profile, evaluated at $10^6$ points per region.}
  \label{sm-tab:euler-results}
  \begin{tabular}{l*{8}{S[table-format=1.2e-1]}}
    \toprule
    & \multicolumn{2}{c}{box} & \multicolumn{2}{c}{$r_c^2+z_c^2<9$}
    & \multicolumn{2}{c}{$r_c=0$} & \multicolumn{2}{c}{$R>R_{\mathrm{ff}}$} \\
    \cmidrule(lr){2-3}\cmidrule(lr){4-5}\cmidrule(lr){6-7}\cmidrule(lr){8-9}
    & {RMSE} & {Max.} & {RMSE} & {Max.} & {RMSE} & {Max.} & {RMSE} & {Max.} \\
    \midrule
    $f_1$ & 7.71e-6 & 1.15e-4 & 9.87e-6 & 4.00e-5 & 1.28e-5 & 4.46e-5 & 1.15e-6 & 1.29e-5 \\
    $f_2$ & 4.65e-5 & 5.73e-4 & 8.37e-6 & 3.83e-5 & 3.30e-5 & 1.14e-4 & 1.15e-6 & 1.90e-5 \\
    $f_3$ & 9.29e-5 & 1.41e-3 & 1.52e-5 & 5.05e-5 & 3.33e-6 & 1.93e-5 & 4.33e-7 & 7.16e-6 \\
    \bottomrule
  \end{tabular}
\end{table}

\bibliographystyle{unsrt}
\bibliography{refs}